\documentclass[12pt, leqno]{amsart}
\usepackage{latexsym,amsmath, amscd, amssymb, amsthm, euscript,
mathrsfs}
\usepackage[all]{xy}
\newtheorem{thm}[subsection]{Theorem}
\newtheorem{cor}[subsection]{Corollary}
\newtheorem{lem}[subsection]{Lemma}

\newtheorem{prop}[subsection]{Proposition}

\theoremstyle{definition}

\newtheorem{defn}[subsection]{Definition}

\theoremstyle{definition}

\theoremstyle{definition}

\newtheorem{rem}[subsection]{Remark}

\numberwithin{equation}{subsection}

\newtheorem{pg}[subsection]{}

\newcommand{\M}{\mathbb{M}}
\newcommand{\Q}{\mathbb{Q}}

\newcommand{\N}{\mathbb{N}}

\newcommand{\isom}{\xrightarrow{\sim}}

\newcommand{\A}{\mathcal{A}}
\newcommand{\mc}{\mathcal }

\newcommand{\Z}{\mathbb{Z}}

\newcommand{\Sp}{\text{\rm Spec}}
\newcommand{\mymargin}[1]{}%{\marginpar{\tiny{#1}}}

\newcommand{\Y}{\mc Y}

\newcommand{\mls}{\mathcal}

  \DeclareMathSymbol{A}{\mathalpha}{operators}{`A}%
  \DeclareMathSymbol{B}{\mathalpha}{operators}{`B}%
  \DeclareMathSymbol{C}{\mathalpha}{operators}{`C}%
  \DeclareMathSymbol{D}{\mathalpha}{operators}{`D}%
  \DeclareMathSymbol{E}{\mathalpha}{operators}{`E}%
  \DeclareMathSymbol{F}{\mathalpha}{operators}{`F}%
  \DeclareMathSymbol{G}{\mathalpha}{operators}{`G}%
  \DeclareMathSymbol{H}{\mathalpha}{operators}{`H}%
  \DeclareMathSymbol{I}{\mathalpha}{operators}{`I}%
  \DeclareMathSymbol{J}{\mathalpha}{operators}{`J}%
  \DeclareMathSymbol{K}{\mathalpha}{operators}{`K}%
  \DeclareMathSymbol{L}{\mathalpha}{operators}{`L}%
  \DeclareMathSymbol{M}{\mathalpha}{operators}{`M}%
  \DeclareMathSymbol{N}{\mathalpha}{operators}{`N}%
  \DeclareMathSymbol{O}{\mathalpha}{operators}{`O}%
  \DeclareMathSymbol{P}{\mathalpha}{operators}{`P}%
  \DeclareMathSymbol{Q}{\mathalpha}{operators}{`Q}%
  \DeclareMathSymbol{R}{\mathalpha}{operators}{`R}%
  \DeclareMathSymbol{S}{\mathalpha}{operators}{`S}%
  \DeclareMathSymbol{T}{\mathalpha}{operators}{`T}%
  \DeclareMathSymbol{U}{\mathalpha}{operators}{`U}%
  \DeclareMathSymbol{V}{\mathalpha}{operators}{`V}%
  \DeclareMathSymbol{W}{\mathalpha}{operators}{`W}%
  \DeclareMathSymbol{X}{\mathalpha}{operators}{`X}%
  \DeclareMathSymbol{Y}{\mathalpha}{operators}{`Y}%
  \DeclareMathSymbol{Z}{\mathalpha}{operators}{`Z}%

\def\CyrillicGuillemets{\DeclareFontEncoding{OT2}{}{}%
     \DeclareFontSubstitution{OT2}{wncyr}{m}{n}%
     \DeclareTextCommand{\guillemotleft}{OT1}{%
        {\fontencoding{OT2}\fontfamily{wncyr}\selectfont\char60}}%
     \DeclareTextCommand{\guillemotright}{OT1}{%
        {\fontencoding{OT2}\fontfamily{wncyr}\selectfont\char62}}}

\def\LasyGuillemets{%
     \DeclareTextCommand{\guillemotleft}{OT1}{\hbox{%
        \fontencoding{U}\fontfamily{lasy}\selectfont(\kern-0.20em(}}%
     \DeclareTextCommand{\guillemotright}{OT1}{\hbox{%
        \fontencoding{U}\fontfamily{lasy}\selectfont)\kern-0.20em)}}}
  \IfFileExists{ot2wncyr.fd}{\CyrillicGuillemets}{\LasyGuillemets}
   \DeclareTextSymbolDefault{\guillemotleft}{OT1}
   \DeclareTextSymbolDefault{\guillemotright}{OT1}
   \def\guill@spacing{\penalty\@M\hskip.8\fontdimen2\font
                               plus.3\fontdimen3\font
                               minus.8\fontdimen4\font}

\newcommand{\ra}{\rightarrow}
\newcommand{\Hom}{\mathrm{Hom}}
\renewcommand{\le}{\textup{lis-{\'e}t}}
\newcommand{\LE}{{\textup{Lisse-{\'e}t}}}
\newcommand{\ET}{{\textup{\'Etale}}}

\newcommand{\et}{\textup{{\'e}t}}

\newcommand{\X}{{\mc X}}
\newcommand{\NN}{{\mathbf N}}

\newcommand{\T}{{\mc T}}
\newcommand{\m}{{\mathfrak m}}
\newcommand{\D}{{\mc D}}
\newcommand{\DD}{{\mathbf D}}
\newcommand{\comp}{\mathrm{Comp}}
\newcommand{\Ab}{\mathrm{Ab}}

\DeclareMathOperator{\Ker}{\mathrm{ker}}
\DeclareMathOperator{\Otimes}{\stackrel{\mathbf L}{\otimes}}
\DeclareMathOperator{\Rhom}{\mathcal{R}\it{hom}}
\DeclareMathOperator{\Rghom}{\hbox{\boldmath{$\mathcal{R}\mathit{hom}_\Lambda$}}}
\DeclareMathOperator{\RgHom}{\hbox{\boldmath{$\mathrm{Rhom}_\Lambda$}}}

\DeclareMathOperator{\egxt}{\hbox{\boldmath{$\mathcal{E}{\it{xt}}_\Lambda$}}}
\DeclareMathOperator{\Egxt}{\hbox{\boldmath{$\mathrm{Ext}_\Lambda$}}}
\DeclareMathOperator{\Hgom}{\hbox{\boldmath{$\mathrm{Hom}_\Lambda$}}}

\DeclareMathOperator{\RHom}{\mathrm{Rhom}}
\DeclareMathOperator{\ext}{\mathcal{E}\it{xt}}
\DeclareMathOperator{\Ext}{\mathrm{Ext}}

\DeclareMathOperator{\Id}{Id}
\begin{document}

\title{The six operations for sheaves on Artin stacks II: Adic Coefficients}

\author{Yves Laszlo and Martin Olsson}
\address{\'Ecole Polytechnique CMLS UMR 7640 CNRS F-91128 Palaiseau Cedex France}
 \email{laszlo@math.polytechnique.fr}
 \address{University of Texas at Austin
Department of Mathematics 1 University Station C1200 Austin, TX
78712-0257, USA}\email{molsson@math.utexas.edu}

\begin{abstract}In this paper we develop a theory of Grothendieck's six operations
for  adic  constructible sheaves on
Artin stacks continuing  the study of the finite coefficients case
in ~\cite{Las-Ols05}.\end{abstract}

\maketitle

\section{Introduction}

In this paper we continue the study of Grothendieck's six operations for sheaves on Artin stacks begun in \cite{Las-Ols05}.  Our aim in this paper is to extend the theory of finite coefficients of loc. cit. to a theory for adic sheaves. In a subsequent paper \cite{pervers} we will use this theory to study perverse sheaves on Artin stacks.

Throughout we work over a regular noetherian scheme $S$ of dimension $\leq 1$. In what follows, all stacks considered will be algebraic locally of finite type over $S$.

Let $\Lambda $ be a complete discrete valuation ring and for every $n$ let $\Lambda _n$ denote the quotient $\Lambda /\mathfrak{m}^n$ so that $\Lambda = \varprojlim \Lambda _n$.
We then define for any  stack $\mls X$
 a triangulated category $\DD_c(\mls X, \Lambda )$ which we call the
 \emph{derived category of constructible $\Lambda $--modules on $\mls X$}
 (of course as in the classical case this is abusive terminology).
 The category $\DD_c(\mls X, \Lambda )$ is obtained from the derived category of projective systems
$\{F_n\}$ of $\Lambda _n$--modules by localizing along the full subcategory of complexes
whose cohomology sheaves are \emph{AR-null} (see \ref{gen-proj} for the meaning of this).
For a morphism $f:\mls X\rightarrow \mls Y$ of finite type of stacks locally of finite type
over $S$ we then define functors
\begin{equation*}
Rf_*:\DD_c^{+}(\mc X, \Lambda )\rightarrow \DD_c^{+}(\mc Y, \Lambda ), \ \
Rf_!:\DD_c^-(\mc X, \Lambda )\rightarrow \DD_c^-(\mc Y, \Lambda ),
\end{equation*}
\begin{equation*}
Lf^*:\DD_c(\mc Y, \Lambda )\rightarrow \DD_c(\mc X, \Lambda ), \ \ Rf^!:\DD_c(\mc Y, \Lambda )\rightarrow \DD_c(\mc X, \Lambda ),
\end{equation*}
\begin{equation*}
\Rghom:\DD_c^{-}(\mc X, \Lambda )^{\text{op}}\times \DD_c^+(\mc X, \Lambda )\rightarrow \DD_c^+(\mc X, \Lambda ),
\end{equation*}
and
\begin{equation*}
(-){\Otimes}(-):\DD_c^-(\mc X, \Lambda )\times \DD_c^-(\mc X, \Lambda )\rightarrow
\DD_c^-(\mc X, \Lambda )
\end{equation*}
satisfying all the usual adjointness properties that one has in
the theory for schemes and the theory for finite coefficients.

In order to develop this theory we must overcome two basic problems. The first one is the
necessary consideration of unbounded complexes which was already apparent in the finite coefficients case. The second one is the non-exactness of the projective limit functor. It should be noted that important previous work has been done
on the subject, especially in ~\cite{Ber03} and~\cite{Eke90} (see also \cite{Jan} for the adic problems). In particular the construction of the normalization functor (\ref{normalization-def}) used in this paper is due to Ekedahl \cite{Eke90}.  None of these works, however, give an entirely satisfactory solution to the problem since for example cohomology with compact support and the duality theory was not constructed.

%%%%%%%%%%%%%%%%%%%%%%%%%%%%%%%%%%%%%%%%%%%%%%
%%%%%%%%%%%%%%%%%%%%%%%%%%%%%%%%%%%%%%%%%%%%

\section{$R\lim$ for unbounded complexes}

Since we are forced to deal with unbounded complexes (in both directions) when considering the functor $Rf_!$ for Artin stacks, we must first collect some results about the unbounded derived category of projective systems of $\Lambda $--modules.  The key tool is \cite{Las-Ols05}, \S2.

\subsection{Projective systems}\label{gen-proj}Let $(\Lambda,\m)$ be a complete local regular ring and
$\Lambda_n=\Lambda/\m^{n+1}$. We demote by $\Lambda_\bullet$ the
pro-ring $(\Lambda_n)_{n\geq 0}$. At this stage, we could have
take any projective system of rings and $\Lambda$ the projective
limit. Let  $\X/S$ be a  stack (by convention algebraic locally of finite type over $S$). For any topos $\mc T$, we will
denote by ${\mc T}^{\NN}$ the topos of projective systems of ${\mc
T}$. These topos will be ringed by $\Lambda,\Lambda_\bullet$
respectively. We denote by $\pi$ the morphism of ringed topos
$\pi:{\mc T}^{\NN}\ra{\mc T}$ defined by $\pi^{-1}(F)=(F)_n$, the
constant projective system. One checks the formula
$$\pi_*=\varprojlim.$$

Recall that for any $F\in \text{Mod}(\T,\Lambda_\bullet)$, the
sheaf $R^i\pi_*F$ is the sheaf associated to  the presheaf
$U\mapsto H^i(\pi^*U,F)$. We'll use several times the fundamental
exact sequence \cite[0.4.6]{Eke85}
\begin{equation}\label{ex-lim-proj}
    0\ra\varprojlim\nolimits^1H^{i-1}(U,F_n)\ra
H^i(\pi^*U,F)\ra\varprojlim H^i(U,F_n)\ra 0.
\end{equation}
 If $*$ denotes the punctual topos, then this sequence is obtained from the Leray spectral sequence associated to the composite
\begin{equation*}
T^{\mathbb{N}}\rightarrow *^{\mathbb{N}}\rightarrow *
\end{equation*}
and the fact that $R^i\varprojlim $ is the zero functor for $i>1$.

Recall that lisse-{\'e}tale topos can be
defined using the lisse-{\'e}tale site $\LE(\X)$ whose objects are
smooth morphisms $U\ra\X$ such that $U$ is an algebraic space of
finite type over $S$.

Recall (cf.~\cite[exp. V]{SGA5}). that a projective system
$M_n,n\geq 0$ in an additive category is \emph{AR-null} if there
exists an integer $r$ such that for every $n$ the composite
$M_{n+r}\ra M_n$ is zero.
\begin{defn} A complex $M$ of $\textup{Mod}({\mc X_\le^\NN},\Lambda_\bullet)$ is\begin{itemize}
    \item \emph{AR-null} if all the ${\mc H}^i(M)$'s are AR-null.
    \item \emph{constructible} if all the ${\mc H}^i(M_n)$'s ($i\in\Z, n\in\NN$) are constructible.
    \item \emph{almost zero} if for any $U\ra \mc X$ in $\LE(\X)$,
the restriction of ${\mc H}^i(M)$ to $\ET(U)$ is AR-null.
\end{itemize}
\end{defn}

Observe that the cohomology sheaves $\mc H^i(M_n)$ of a
constructible complex are by definition cartesian.

\begin{rem} A constructible complex $M$ is almost zero if and only
if its restriction to some presentation $X\ra\X$ is almost zero,
meaning that there exists a covering of $X$ by open subschemes $U$
of finite type over $S$ such that the restriction $M_U$ of $M$ to
$\ET(U)$ is AR-null.\end{rem}

\subsection{Restriction of $R\pi_*$ to $U$} Let $U\ra\X$ in
$\LE(\X)$. The restriction of a complex $M$ of $\X$ to $\ET(U)$ is
denoted as usual $M_U$.
\begin{lem}\label{commutU} One has $R\pi_*(M_U)=(R\pi_*M)_U$ in
$\D(U_\et,\Lambda)$.
\end{lem}
\begin{proof} We view $U$ both as a sheaf on $\X$ or as the constant
projective system $\pi^*U$. With this identification, one has
$(\X_{\le|U})^\NN=(\X_{\le}^\NN)_{|U}$ which we will denote by
$\X_{\le|U}^\NN$. The following diagram commutes
$$\xymatrix{\X_{\le|U}^\NN\ar[r]^j\ar[d]^\pi&\X_\le^\NN\ar[d]^\pi\\\X_{\le|U}\ar[r]^j&\X_\le}$$
where $j$ denotes the localization morphisms and $\pi$ is as
above. Because the left adjoint $j_!$ of $j^*$ is exact, $j^*$
preserves K-injectivity. We get therefore
\begin{equation}\label{eqj}
    R\pi_*j^*=j^*R\pi_*.
\end{equation}
 As
before, the morphism of sites
$\epsilon^{-1}:\ET(U)\hookrightarrow\LE(\X)_{|U}$ and the
corresponding one of total sites induces a commutative diagram
$$\xymatrix{\X_{\le|U}^\NN\ar[r]^\epsilon\ar[d]^\pi&U_\et^\NN\ar[d]^\pi\\\X_{\le|U}\ar[r]^\epsilon&U_\et}$$

Since $\epsilon _*$ is exact with an exact left adjoint $\epsilon ^*$, one has
\begin{equation}\label{eqeps}
    R\pi_*\epsilon_*=\epsilon_*R\pi_*.
\end{equation}

One gets therefore \begin{eqnarray*}
  (R\pi_*M)_U &=& \epsilon_*j^*R\pi_*M \\
   &=& \epsilon_*R\pi_*j^*M \textup{ by }(\ref{eqj}) \\
   &=& R\pi_*\epsilon_*j^*M \textup{ by }(\ref{eqeps}) \\
   &=& R\pi_*(M_U)
\end{eqnarray*}

\end{proof}

As before, let $\T$ be a topos and
let $\A$ denote the category of $\Lambda _\bullet $--modules in $\T^{\mathbb{N}}$.

\begin{prop}[Lemma 1.1 of Ekedahl]\label{E1} Let $M$ be a complex of $\A$.
\begin{enumerate}
    \item If $M$ is AR-null, then
${R}\pi_*M=0$.
    \item If $M$ almost zero, then
${R}\pi_*M=0$.
\end{enumerate}    \end{prop}
\begin{proof} Assume $M$ is AR-null. By \cite[lemma 1.1]{Eke90}
$R\pi_*\mc H^j(M)=0$ for all $j$. By \cite[2.1.10]{Las-Ols05} one
gets ${R}\pi_*M=0$. The second point follows from (1) using
\ref{commutU}.\end{proof}

\begin{lem}[Lemma 1.3 iv) of Ekedahl]\label{E2} Let $M$ be  complex in ${\mc T}$ of $
\Lambda_n$-modules. Then,  the adjunction morphism
$M\ra{{R}}\pi_*\pi^*M$ is an isomorphism.
\end{lem}
\begin{rem} Here we view $\pi $ is a morphism of ringed topos $(\mc T^\NN, \Lambda _n)\rightarrow (\T, \Lambda _n)$.  Then functor $\pi ^*$ sends a $\Lambda _n$--module $M$ to the constant projective system $M$.  In particular, $\pi ^*$ is exact (in fact equal to $\pi ^{-1}$) and hence passes to the derived category.
\end{rem}
\noindent {\it Proof of \ref{E2}:}
The sheaf $R^i\pi _*{\mc H}^j(\pi ^*M)$ is the sheaf associated to the presheaf sending $U$ to $H^i(\pi ^*U,  \mc H^j(\pi ^*M))$.  It follows from \ref{ex-lim-proj} and the fact that the system $H^{i-1}(U, \mc H^j(\pi ^*M)_n)$ satisfies the Mittag-Leffler condition that this presheaf is isomorphic to the sheaf associated to the presheaf $U\mapsto \varprojlim H^i(U, \mc H^j(\pi ^*M)_n) = H^i(U, \mc H^j(\pi ^*M)$.  It follows that
 ${{R}}^i\pi_*{\mc H}^j(\pi^*M)=0$ for all $i>0$
 and  $$\mc H^jM=R\pi_*{\mc
H}^j(\pi^*M).\leqno(*)$$ By \cite[2.1.10]{Las-Ols05} one can
therefore assume $M$ bounded from below. The lemma follows
therefore by induction from (*) and from the distinguished
triangles
$$\mc H^{j}[-j]\ra\tau_{\geq j}M\ra\tau_{\geq j+1}M.$$\qed

In fact, we have the following stronger result:
\begin{prop}\label{2.9} Let $N\in \D(\T^\NN, \Lambda _n)$ be a complex of projective systems such that for every $m$ the map
\begin{equation}\label{transmaps}
N_{m+1}\rightarrow N_m
\end{equation}
is a quasi--isomorphism.  Then the natural map $\pi ^*R\pi _*N\rightarrow N$ is an isomorphism. Consequently, the functors $(\pi ^*, R\pi _*)$ induce an equivalence of categories between $\D(\T, \Lambda _n)$ and the category of complexes $N\in \D(\T^\NN, \Lambda _n)$ such that the maps \ref{transmaps} are all isomorphism.
\end{prop}
\begin{proof} By \cite[2.1.10]{Las-Ols05} it suffices to prove that the map $\pi ^*R\pi _*N\rightarrow N$ is an isomorphism for $N$ bounded below.  By devissage using the distinguished triangles
$$
\mc H^j(N)[j]\rightarrow \tau _{\geq j}N\rightarrow \tau _{\geq j+1}N
$$
one further reduces to the case when $N$ is a constant projective system of sheaves where the result is standard (and also follows from \ref{E2}.
\end{proof}

\section{$\lambda$-complexes} Following Behrend and~\cite[exp. V, VI]{SGA5}, let us start with a definition.
Let $\mc X$ be an algebraic stack locally of finite type over $S$, and let $\A$ denote the category of $\Lambda _\bullet $--modules in $\mc X_{\le}^{\mathbb{N}}$.

\begin{defn}  We say that\begin{itemize}
    \item a system $M=(M_n)_n$ of $\A$
    is \emph{adic} if all the $M_n$'s are
    constructible and moreover all morphisms
    $$\Lambda_n\otimes_{\Lambda_{n+1}}M_{n+1}\ra M_n$$ are
    isomorphisms; it is called  \emph{almost adic}  if all the $M_n$'s are
    constructible and if for every $U$ in $\LE(\X)$ there is a morphism
    $N_U\ra M_U$ with almost zero kernel and cokernel with $N_U$ adic in $U_\et$.
    \item a complex  $M=(M_n)_n$ of $\A$ is called a \emph{$\lambda$-complex}  if all the cohomology
    modules ${\mc H}^i(M)$ are almost adic. Let $\D_c(\A)\subset \D(\A)$ denote the full triangulated subcategory whose objects are $\lambda $--complexes. The full subcategory of $\D_c(\A)$ of complexes concentrated in degree $0$ is called the category of
    $\lambda$-\emph{modules}.
    \item The category $\DD_c(\X,\Lambda)$ (sometimes written just $\DD_c(\X)$ if the reference to $\Lambda $ is clear) is the quotient of the
category $\D_c(\A)$  by the full subcategory of almost zero complexes.
    \end{itemize}
\end{defn}

\begin{rem} Let $X$ be a noetherian scheme.  The condition that a sheaf of $\Lambda _\bullet $--modules $M$  in $X_\et ^\mathbb{N}$ admits a morphism $N\rightarrow M$ with $N$ adic is \'etale local on $X$.  This follows from \cite[V.3.2.3]{SGA5}. Furthermore, the category of almost adic $\Lambda _\bullet $-modules is an abelian subcategory closed under extensions (a Serre subcategory) of the category of all $\Lambda _\bullet $--modules in $X_\et ^{\mathbb{N}}$.  From this it follows that for an algebraic stack $\X$, the category of almost adic $\Lambda _\bullet $--modules is a Serre subcategory of the category of all $\Lambda _\bullet $--modules in $\mls X_{\le }^{\mathbb{N}}$.

In fact if $M$ is almost adic on $X$, then the pair $(N, u)$ of an adic sheaf $N$ and an AR--isomorphism $u:N\rightarrow M$ is unique up to unique isomorphism.  This follows from the description in \cite[V.2.4.2 (ii)]{SGA5} of morphisms in the localization of the category of almost adic modules by the subcategory of AR-null modules.  It follows that even when $X$ is not quasi--compact, an almost adic sheaf $M$ admits a morphism $N\rightarrow M$ with $N$ adic whose kernel and cokernel are AR--null when restricted to an quasi--compact \'etale $X$--scheme.
\end{rem}

As usual, we denote by $\Lambda$ the image of
$\Lambda_\bullet$ in $\DD_c(\X)$. By~\cite[exp {V}]{SGA5}
the quotient of the subcategory of almost adic modules by the category
of almost zero modules is abelian. By construction, a morphism
$M\ra N$ of $\D_c(\A)$ is an isomorphism in
$\DD_c(\X)$ if and only if its cone is almost zero. $\DD_c(\X)$ is a
triangulated category and has a natural $t$-structure whose heart
is the localization of the category of $\lambda$-modules by the full subcategory of almost zero systems (cf.~\cite{Ber03}). Notice
however that we do not know at this stage that in general
$\Hom_{\DD_c(\X)}(M,N)$ is a (small) set. In fact, this is equivalent
to finding a left adjoint of the projection $\D_c(\A)\ra\DD_c(\X)$
\cite[section 7]{Nee01}. Therefore, we have to find a
normalization functor $M\ra\hat M$. We'll prove next that a
suitably generalized version of Ekedahl's functor defined in \cite{Eke90} does the job.
Note that by \ref{E1} the functor $R\pi _*:\D_c(\A)\rightarrow \D_c(\X)$ factors uniquely through a functor which we denote by the same symbols $R\pi _*:\DD_c(\X)\rightarrow \D_c(\X)$.

\begin{defn}\label{normalization-def} We define the \emph{normalization functor}
\begin{equation*}
\DD_c(\X)\rightarrow \D(\A), \ \ M\mapsto \hat M
\end{equation*}
by the formula $\hat
M=L\pi^*{R}\pi_*M$. A complex $M\in \D(\A)$ is \emph{normalized} if the natural map $\hat M\rightarrow M$ is an isomorphism (where we write $\hat M$ for the normalization functor applied to the image of $M$ in $\DD_c(\X)$).\end{defn}
 Notice that $\hat\Lambda=\Lambda$ (write $\Lambda_\bullet=L\pi^*\Lambda$ and
 use~\ref{normalchar} below for instance).

\begin{rem} Because $\Lambda$ is regular, the Tor dimension of
$\pi$ is $d=\dim(\Lambda)<\infty$ and therefore we do not have to
use Spaltenstein's theory in order to define $\hat M$.\end{rem}

\begin{prop}[{\cite[2.2 (ii)]{Eke90}}]\label{normalchar} A complex $M\in \D(\X_\le^\NN, \Lambda _\bullet )$ is normalized if and only if for all $n$ the natural map
\begin{equation}\label{transition}
\Lambda _n\Otimes _{\Lambda _{n+1}}M_{n+1}\rightarrow M_n
\end{equation}
is an isomorphism.
\end{prop}
\begin{proof}
If $M = L\pi ^*N$ for some $N\in \D(\X_\le, \Lambda )$ then for all $n$ we have $M_n = \Lambda _n\Otimes _\Lambda N$ so in this case the morphism \ref{transition} is equal to the natural isomorphism
$$
\Lambda _n\Otimes _{\Lambda _{n+1}}\Lambda _{n+1}\Otimes _\Lambda N\rightarrow \Lambda _n\Otimes _\Lambda N.
$$
This proves the ``only if'' direction.

For the ``if'' direction, note that since the functors $e_n^*$ form a conservative set of functors, to verify that $\hat M\rightarrow M$ is an isomorphism it suffices to show that for every $n$ the map $e_n^*\hat M\rightarrow e_n^*M$ is an isomorphism. Equivalently we must show that the natural map
$$
\Lambda _n\Otimes _\Lambda R\pi _*(M)\rightarrow M_n
$$
is an isomorphism.  As discussed in \cite[bottom of page 198]{Eke90}, the natural map $L\pi ^*\Lambda _n\rightarrow \pi ^*\Lambda _n$ as AR--null cone.  In the case when $\Lambda $ is a discrete valuation ring with uniformizer $\lambda $, this can be seen as follows. A projective resolution of $\Lambda _n$ is given by the complex
$$
\begin{CD}
\Lambda @>\times \lambda ^{n+1}>> \Lambda.
\end{CD}
$$
>From this it follows that $L\pi ^*(\Lambda _n)$ is represented by the complex
$$
\begin{CD}
(\Lambda _m)_m@>\times \lambda ^{n+1}>> (\Lambda _m)_m.
\end{CD}
$$ Therefore the cone of $L\pi ^*(\Lambda _n)\rightarrow \pi ^*\Lambda _n$ is in degrees $m\geq n$ and up to a shift equal to $\lambda ^{m-n}\Lambda _m$ which is  AR-null.

Returning to the case of general $\Lambda $, we obtain from the  projection formula and \ref{E1}
$$
\Lambda _n\Otimes _\Lambda R\pi _*(M)\simeq e_n^{-1}R\pi _*(L\pi ^*\Lambda _n\Otimes M)\simeq e_n^{-1}R\pi _*(\pi ^*\Lambda _n\Otimes _{\Lambda _\bullet }M) = R\pi _*(\Lambda _n\Otimes _{\Lambda _\bullet }M).
$$
The proposition then follows from \ref{2.9}.
\end{proof}

 We have a localization result analogous
to lemma~\ref{commutU}. Let $M\in \D(\X_\le,\Lambda)$.

\begin{lem}\label{commutLU} One has $L\pi^*(M_U)=(L\pi^*M)_U$ in
$\D(U_\et^\NN,\Lambda_\bullet)$.
\end{lem}

\begin{proof} We use the notations of the proof of lemma~\ref{commutU}.
First,  $j^*=Lj^*$ commutes with $L\pi^*$ due to the commutative
diagram
$$\xymatrix{\X_{\le|U}^\NN\ar[r]^j\ar[d]^\pi&\X_\le^\NN\ar[d]^\pi\\\X_{\le|U}\ar[r]^j&\X_\le}.$$

One is therefore reduced to prove that $\epsilon_*=R\epsilon_*$
commutes with $L\pi^*$. We have certainly, with a slight abuse of
notation,
$$\epsilon^{-1}\Lambda_\bullet=\Lambda_\bullet.$$
Therefore, if $N$ denotes the restriction of $M$ to $\X_{|U}$ we get

\begin{eqnarray*}
  \epsilon_*L\pi^*N &=& \epsilon_*(\Lambda_\bullet\Otimes_{\pi^{-1}\Lambda}\pi^{-1}N)\\
   &=& \epsilon_*(\epsilon^{-1}\Lambda_\bullet\Otimes_{\pi^{-1}\Lambda}\pi^{-1}N)\\
   &=& \Lambda_\bullet\Otimes_{\epsilon_*\pi^{-1}\Lambda}\epsilon_*\pi^{-1}N
   \textup{  by the projection formula}\\
   &=& \Lambda_\bullet\Otimes_{\pi^{-1}\Lambda}\pi^{-1}\epsilon_*N
   \textup{ because $\epsilon_*$ commutes with $\pi^{-1}$} \\
   &=& L\pi^*\epsilon_*N.
\end{eqnarray*}
\end{proof}

\begin{rem}\label{hatMbounded} The same arguments used in the proof shows that if
$M\in\D_c(\X_\le,\Lambda_\bullet)$ and $M_U$ is bounded for $U\in \LE(\X)$, then
$\hat M_U$ is also bounded. In particular, all $\hat
M_{U,n}$ are of finite tor-dimension.
\end{rem}

 \begin{cor}\label{chapeauM-M} Let  $M\in\D(\X_\le^\NN,\Lambda_\bullet)$ and $U\ra\X$ in
$\LE(\X)$.Then, the adjunction morphism
$$\hat M\ra M$$ restricts on $U_\et$ to the adjunction morphism
$L\pi^*R\pi_{ *}M_U\ra M_U$.\end{cor} \begin{proof} It is an
immediate consequence of~lemmas~\ref{commutLU}
and~\ref{commutU}.\end{proof}

We assume now that $\Lambda$ is a discrete valuation ring with
uniformizing parameter $\lambda$. Let us prove the analogue
of~\cite[Proposition 2.2]{Eke90}.

%\begin{prop}\label{hat-cart}Let
%$M\in\D_{\mathrm{cart}}(\X_\le^\NN,\Lambda_\bullet)$ be a
%cartesian complex (namely $ \mc H^i(M)$ cartesian for each $i$).
%Then $\hat M$ is cartesian. Moreover, if $\X$ is of finite type,
%the normalizaton preserves $\D_{\mathrm{cart}}^{\pm,b}$.
%\end{prop}
%\begin{proof} Immediate consequence of~\ref{chapeauM-M} and of the
%finite cohomological dimension of $R\pi_{U*}$ and $L\pi_{U*}$ for
%a presentation $U\ra\X$ if $\X$ is of finite type.
%\end{proof}

%%%%%%%%%%%%%%%%%%

\begin{thm}\label{th-normalisation}Let $M$ be a $\lambda$-complex.
Then, $\hat M$ is constructible and $\hat M\ra M$ has an almost
zero cone.
\end{thm}

\begin{proof}  Let $U\ra\X$ in
$\LE(\X)$ and $$N=M_U\in\D_c(U_\et,\Lambda_\bullet).$$ Let us
prove first that $(\hat M)_U\in\D(U_\et)$ is constructible and
that the cone of $(\hat M)_U\ra M_U$ is $AR$-null. We proceed by
successive reductions.
\begin{enumerate}
    \item  Let
$d_U=\text{cd}_{\ell}(U_\et)$ be the $\ell$-cohomological
dimension of $U_\et$. By an argument similar to the one used in the proof of \ref{E2} using \ref{ex-lim-proj}, the cohomological dimension of
$R\pi_*$ is  $\leq 1+d_U$. Therefore, $R\pi_*$
 maps
$\D^{\pm,b}(U_\et^\NN)$ to $\D^{\pm,b}(U_\et)$. Because $L\pi^*$ is of finite
cohomological dimension, the same is true for the normalization
functor. More preciseley, there exists an integer $d$ (depending
only on $U\ra\X$ and $\Lambda$) such that for every $a$

$$N\in
    \D^{\geq a}(U_\et ^\NN)\Rightarrow\hat N\in \D^{\geq a-d}(U_\et^\NN)\textup{ and }N\in
    \D^{\leq a}(U_\et ^\NN)\Rightarrow\hat N\in \D^{\leq a+d}(U_\et^\NN).$$
    \item One can assume $N\in \text{Mod}(U_\et,\Lambda_\bullet)$.
    Indeed, one has by the previous observations $$\mc H^i(\hat
    N)=\mc H^i(\widehat{N_i})$$
    where $N_i=\tau_{\geq i-d}\tau_{\leq i+d}N$. Therefore one can assume $N$ bounded.
    By induction, one can
    assume $N$ is a $\lambda$-module.
    \item One can assume $N$ adic. Indeed, there exists a morphism $A\ra N$
    with AR-null kernel and cokernel with $A$ adic. In particular the cone of $A\ra N$ is AR-null.
    It is therefore enough to observe that $\hat A=\hat N$, which is a consequence of~\ref{E1}.
    \item We use without further comments basic
facts about the abelian category of $\lambda$-modules
(cf.~\cite[exp. V]{SGA5} and \cite[Rapport sur la formule des traces]{SGA4.5}).

In the category of $\lambda$-modules, there exists $n_0$ such that
$N/\Ker(\lambda^{n_0})$ is torsion free (namely the action of
$\lambda$ has no kernel). Because $\DD_c(U_\et)$ is triangulated, we just
have to prove that the normalization of both
$N/\Ker(\lambda^{n_0})$ and $\Ker(\lambda^{n_0})$ are
constructible and the corresponding cone is AR-null.

\item The case of ${\bar N}=N/\Ker(\lambda^{n_0})$. An  adic
representative $L$ of ${\bar N}$ has flat components $L_n$, in
other words
$$\Lambda_n\Otimes_{\Lambda_{n+1}}L_{n+1}\ra L_n$$ is an
isomorphism. By~\ref{normalchar}, $L$ is
normalized and therefore $\widehat{{\bar N}}=\hat L=L$
is constructible (even adic) and the cone $L=\hat {\bar N}\ra {\bar
N}$ is AR-null because the kernel and cokernel of $L\ra \bar N$
are AR-null.

\item

We can therefore assume $\lambda^{n_0}N=0$ (in the categories of
$\lambda$-modules up to AR-isomorphisms) and even $\lambda N=0$
(look at the $\lambda$-adic filtration). The morphism
\begin{equation}\label{N-N0}
    (N_n)_{n\in\NN}\ra (N_n/\lambda N_n)_{n\in\NN}
\end{equation} has
AR-zero kernel and the normalization of both are therefore the
same. But, $N$ being adic, one has $N_n/\lambda N_n=N_{0}$ for
$n\geq 0$. In particular, the morphism~\ref{N-N0} is nothing but
\begin{equation}\label{N-N1}
    N\ra\pi^*N_0
\end{equation} and is an AR-isomorphism  and
$$\hat N=\widehat{\pi^*N_0}=L\pi^*N_{0}$$  (\ref{E2}). One
therefore has to show that the cone $C$ of $L\pi^*N_0\ra \pi^*N_0$
is almost zero. As before, one can assume replace $\X_\le$ by
$U_\et$ for some affine scheme $U$. On $U$, there exists a
\emph{finite} stratification on which $N_0$ is smooth. Therefore,
one can even assume that $N_0$ is constant and finally equal to
$\Lambda_0$. In this case the cone of $L\pi ^*\Lambda _0\rightarrow \pi ^*\Lambda _0$ is AR-null by the same argument used in the proof of \ref{normalchar}, and therefore this proves the first point.
\end{enumerate}
We now have to prove that $\hat M$ is cartesian.
By~\ref{chapeauM-M} again, one is reduced to the following
statement :

Let $f:V\ra U$ be an $\X$-morphism in $\LE(\X)$ which is smooth.
Then, $$f^*\widehat{ M_U}=\widehat{ M_V}=\widehat
{f^*M_U}\footnote{By~\ref{chapeauM-M}, the notations there is no
ambiguity in the notation.}.$$

The same reductions as above allows to assume that $M_U$ is
concentrated in degree $0$, and that we have  a distinguished
triangle $$L\ra M_U\ra C$$ with $C$ AR-null and $L$ either equal
to $\Lambda_0$ or adic with flat components. Using the exactness
of $f^*$ and the fact that $M$ is cartesian, one gets a
distinguished triangle
$$f^*L\ra M_V\ra f^*C$$ with $f^*C$ AR-null. We get therefore
$f^*\hat M_U=f^*\hat L$ and $\widehat { f^*M_U}=\widehat {f^* L}$ : one can
assume $M_U=L$ and $M_V=f^*L$. In both cases, namely $L$ adic with
flat components or $L=\Lambda_0$, the computations above shows
$f^*\hat L=\widehat{ f^* L}$ proving that $\hat M$ is cartesian.
\end{proof}

\begin{rem} The last part of the proof of the first point is proved in a greater
generality in~\cite[lemma 3.2]{Eke90}.  \end{rem}

\begin{rem}
In general the functor $R\pi _*$ does not take cartesian sheaves to cartesian sheaves.  An example suggested by J. Riou is the following: Let $Y = \Sp (k)$ be the spectrum of an algebraically closed field and $f:X\rightarrow Y$ a smooth $k$--variety.  Let $\ell $ be a prime invertible in $k$ and let $M = (M_n)$ be the projective system $\Z/\ell^{n+1}$ on $Y$.  Then $R\pi _*M$ is the constant sheaf $\Z_\ell $, and so the claim $R^i\Gamma (f^*R\pi _*M)$ is the cohomology of $X$ with values in the constant sheaf $\Z_\ell $.  On the other hand,  $R^i\Gamma (R\pi _*(f^*M))$ is the usual $\ell $--adic cohomology of $X$ which in general does not agree with the cohomology with coefficients in $\Z_\ell $.
\end{rem}

\begin{cor}\label{calcul-pou-rhom} Let $M\in \DD_c(\X)$. Then for any $n\geq 0$,
one has
$$\Lambda_n\Otimes_\Lambda R\pi_*M\in\D_c(\X,\Lambda_n).$$\end{cor}
\begin{proof} Indeed, one has $e_n^{-1}\hat M=\Lambda_n\Otimes_\Lambda
R\pi_*M$ which is constructible
by~\ref{th-normalisation}.\end{proof}

We are now able to prove the existence of our adjoint.

\begin{prop}\label{small-hom} The normalization functor is a left adjoint of
the projection $\D_c(\X^\NN)\ra\DD_c(\X)$. In particular,
$\Hom_{\DD_c(\X)}(M,N)$ is small for any $M,N\in\DD_c(\X)$.
\end{prop}
\begin{proof} With a slight abuse of notations, this means
$\Hom_{\D_c}(\hat M,N)=\Hom_{\DD_c(\X)}(M,N)$. If we start with a
morphism $\hat M\ra N$, we get a diagram $$\xymatrix{&\hat
M\ar[ld]\ar[rd]\\M&&N}$$ where $\hat M\ra M$ is an isomorphism in
$\DD_c(\X)$ by~\ref{chapeauM-M} and \ref{th-normalisation} which
defines a morphism in $\Hom_{\DD_c(\X)}(M,N)$. Conversely, starting
from a diagram $$\xymatrix{&L\ar[ld]\ar[rd]\\M&&N}$$ where $L\ra
M$ is an isomorphism in $\DD_c(\X)$.  Therefore one has  $\hat M=\hat L$
(\ref{E1}), and  we get a morphism $\hat M\ra\hat N$ in
$\D_c$  and therefore, by composition, a morphism $\hat M\ra N$.
One checks that these construction are inverse each
other.\end{proof}

%%%%%%%%%%%%%%%%%%%%%%%

\subsection{Comparison with Deligne's approach} Let
$M,N\in\D_c(\X^\NN, \Lambda _\bullet)$ and assume $M$ is normalized. Then there is a
sequence of morphisms
$$\RHom(M_n,N_n)\ra\RHom(M_n,N_{n-1})=\RHom(\Lambda_{n-1}\Otimes_{\Lambda_n}M_n,N_{n-1})=
\RHom(M_{n-1},N_{n-1}).$$ Therefore, we get for each $i$ a
projective system $(\Ext^i(M_n,N_n))_{n\geq 0}$.
\begin{prop}\label{compa-deligne} Let $M,N\in\D_c(\X^\NN)$ and assume $M$ is normalized.
Then there is an exact sequence $$
0\ra\varprojlim\nolimits^1\Ext^{-1}(M_n,N_n)\ra\Hom_{\D_c(\X)}(M,N)\ra\varprojlim
\Hom_{\D(\X,\Lambda_n)}(M_n,N_n)\ra 0.$$
\end{prop}

\begin{proof} Let $\X_\le^{\leq n}$ be the $[0\cdots n]$-simplicial topos
of projective systems $(F_m)_{m\leq n}$ on $\X_\le$. Notice that
the inclusion $[0\cdots n]\ra\NN$ induces an open immersion of the
corresponding topos and accordingly an open immersion
\begin{equation}\label{defjn}j_n:\X_\le^{\leq n}\hookrightarrow \X_\le^\NN.\end{equation}
The inverse image functor is just the truncation $F=(F_m)_{m\leq
0}\mapsto F^{\leq n}=(F_m)_{m\leq n}$. We get therefore an
inductive system of open sub-topos of $\X_\le^\NN$ :
$$\X_\le^{\leq 0}\hookrightarrow\X_\le^{\leq
1}\hookrightarrow\cdots\X_\le^{\leq
n}\hookrightarrow\cdots\X_\le^\NN.$$  Fixing $M$,
let
$$
F:\comp^+(\X^\NN,\Lambda_\bullet)\rightarrow \Ab^\NN
$$
be the functor
$$
N\mapsto \Hom(M^{\leq
n},N^{\leq n}))_n.
$$
Then there is  a
commutative diagram
$$\xymatrix{
\comp^+(\X^\NN,\Lambda_\bullet)\ar[rrr]^-{F}
\ar[rrd]_{\Hom(M,-)} &&&\Ab^\NN\ar[ld]^{\varprojlim}\\
&&\Ab&}$$ which yields the equality $$\RHom(M,N)=R\varprojlim\circ
RF(N).$$ By the definition of $F$ we have $R^qF(N)=\Ext^q(M^{\leq n},N^{\leq n}))_n$.
Because $\varprojlim$ is of cohomological dimension $1$, there is an equality of functors
$\tau_{\geq 0}R\varprojlim=\tau _{\geq 0}R\varprojlim\tau_{\geq -1}$.

Using the distinguished triangles $$(\mc
H^{-d}RF(N))[d]\ra\tau_{\geq -d} RF(N)\ra\tau_{\geq -d+1}RF(N)$$
we get for $d=1$ an exact sequence
$$0\ra\lim\nolimits^1\Ext^{-1}(M^{\leq n},N^{\leq n})\ra\Hom(M,N)\ra
R^0\varprojlim\tau_{\geq 0}RF(N)\ra 0,$$  and for
$d=0$
$$\varprojlim\Hom(M^{\leq n},N^{\leq n})=R^0\varprojlim\tau_{\geq
0}RF(N)^.$$ Therefore we just have to show the formula $$\Ext^q(M^{\leq
n},N^{\leq n})=\Ext^q(M_n,N_n)$$
which follows from the following lemma which will also be useful below.
\end{proof}

\begin{lem}\label{calc-loc} Let $M,N\in\D(\X^\NN)$ and assume $M$ is normalized.
Then, one has\begin{enumerate}
    \item $\RHom(M^{\leq n},N^{\leq n})=\RHom(M_n,N_n)$.
    \item $e_n^{-1}\Rhom(M,N)=\Rhom(M_n,N_n)$.
\end{enumerate}
\end{lem}
\begin{proof}
Let $\pi_n:\X_\le^{\leq n}\ra\X_\le$ the restriction of $\pi$. It
is a morphism of ringed topos ($\X$ is ringed by $\Lambda_n$ and
$\X_\le^{\leq n}$ by $j_n^{-1}(\Lambda_\bullet)=(\Lambda_m)_{m\leq
n}$). The morphisms $e_i:\X\ra\X^\NN, i\leq n$ can be localized in
    $\tilde e_i:\X\ra\X^{\leq n}$, characterized by $e_i^{-1}M^{\leq n}=M_i$
    for any object $M^{\leq n}$ of $\X_\le^{\leq n}$.They form a conservative
    sets of functors satisfying \begin{equation}\label{ej}
    e_i=j_n\circ\tilde
    e_i
\end{equation} One has
$$\pi_{n*}(M^{\leq n})=\varprojlim_{m\leq n}M_m=M_n=\tilde e_n^{-1}(M).$$ It follows that
$\pi_{n*}$ is exact and therefore

\begin{equation}\label{Rpi=e}
R\pi_{n*}=\pi_{n*}=\tilde e_n^{-1}.
\end{equation}
The isomorphism $M_n\ra R\pi_{n*}M^{\leq n}$ defines by adjunction
a morphism $L\pi^*M_n\ra M^{\leq n}$ whose pull back by $\tilde
e_i$ is $\Lambda_i\Otimes_{\Lambda_n}M_n\ra M_i$. Therefore, one
gets\begin{equation}\label{Lpi=M}
    L\pi_n^*M_n=M^{\leq n}
\end{equation}
because $M$ is normalized. Let us prove the first point. One has
$$\RHom(M^{\leq n},N^{\leq n})\stackrel{\ref{Lpi=M}}{=}\RHom(L\pi_{n}^*M_n,N^{\leq
n})\stackrel{\text{adjunction}}{=}\RHom(M_n,R\pi_{n*}N^{\leq
n})\stackrel{\ref{Rpi=e}}{=}\RHom(M_n,N_n).$$
 proving the first point. The second point is analogous :
\begin{eqnarray*}
  e_n^{-1}\Rhom(M,N) &=&\tilde e_n^{-1}j_n^{-1}\Rhom(M,N) \ (\ref{ej}) \\
   &=& \tilde e_n^{-1}\Rhom(M^{\leq n},N^{\leq n}) \ (j_n\text{ open immersion})\\
   &=&R\pi_{n*}\Rhom(M^{\leq n},N^{\leq n}) \ (\ref{Rpi=e}) \\
   &=& R\pi_{n*}\Rhom(L\pi_n^*M_n,N^{\leq n}) \ (\ref{Lpi=M})  \\
  &=& \Rhom(M_n,R\pi_{n*}N^{\leq n})\ (\text{projection formula)} \\
     &=& \Rhom(M_n,N_n)\ (\ref{Rpi=e})\\
\end{eqnarray*}
\end{proof}

\begin{cor}\label{cor-compa-deligne} Let $M,N\in\D_c(\X^\NN)$ be normalized complexes.
Then, one has an exact sequence $$
0\ra\varprojlim\nolimits^1\Ext^{-1}(M_n,N_n)\ra\Hom_{\DD_c(\X)}(M,N)\ra\varprojlim
\Hom_{\D(\X,\Lambda_n)}(M_n,N_n)\ra 0.$$
\end{cor}

\begin{rem}\label{rem-cont-coh} Using similar arguments (more precisely
using Grothendieck spectral sequence of composite functors rather
than truncations as above), one can show  that for any adic
constructible sheaf $N$, the cohomology group $H^*(\mc
X,N)\stackrel{\text{def}}{=}\Ext^*_{\DD_c(\X)}(\Lambda,N)$
coincides with the continuous cohomology group of~\cite{Jan}
(defined as the derived functor of $N\mapsto\varprojlim H^0(\mc
X,N_n)$).
\end{rem}

 Now
let $k$ be either a finite field or an algebraically closed field,
set $S = \Sp (k)$, and let $X$ be a $k$--variety. In this case
Deligne defined in \cite[1.1.2]{WeilII} another triangulated
category which we shall denote by $\D_{c, \text{Del}}^b(X, \Lambda
)$. This triangulated category is defined as follows.  First let
$\D^-_{\text{Del}}(X, \Lambda )$ be the 2--categorial projective
limit of the categories $\D^-(X, \Lambda _n)$ with respect to the
transition morphisms
$$
\Otimes _{\Lambda _n}\Lambda _{n-1}:\D^-(X, \Lambda _n)\rightarrow \D^-(X, \Lambda _{n-1}).
$$
So an object $K$ of $\D^-_{\text{Del}}(X, \Lambda )$ is a projective system $(K_n)_n$ with each $K_n\in D^-(X, \Lambda _n)$ and isomorphisms $K_n\Otimes _{\Lambda _n}\Lambda _{n-1}\rightarrow K_{n-1}$.  The category $\D_{c, \text{Del}}^b(X, \Lambda )$ is defined to be the full subcategory of $\D^-_{\text{Del}}(X, \Lambda )$ consisting of objects $K = (K_n)$ with each $K_n\in \D_c^b(X, \Lambda _n)$. By \cite[1.1.2 (e)]{WeilII} the category $\D^b_{c, \text{Del}}(X, \Lambda )$ is triangulated with distinguished triangles defined to be those triangles inducing distinguished triangles in each $\D^b_c(X, \Lambda _n)$.

By \ref{normalchar}, there is a natural triangluated functor
\begin{equation}\label{comparisonfunctor}
F:\DD_c(X, \Lambda )\rightarrow \D_{c, \text{Del}}^b(X, \Lambda ), \ \ M\mapsto \hat M.
\end{equation}

\begin{lem} Let $K = (K_n)\in \D^b_{c, \text{\rm Del}}(X, \Lambda )$ be an object.

(i) For any integer $i$, the projective system $(\mc H^i(K_n))_n$ is almost adic.

(ii) If $K_0\in \D^{[a, b]}_c(X, \Lambda _0)$, then for $i<a$ the system $(\mc H^i(K_n))_n$ is AR-null.
\end{lem}
\begin{proof} By the same argument used in \cite[1.1.2 (a)]{WeilII} it suffices to consider the case when $X= \Sp (k)$.  In this case, there exists by \cite[XV, page 474 Lemme 1 and following remark]{SGA5} a bounded above complex of finite type flat $\Lambda $--modules $P^\cdot $ such that $K$ is the system obtained from the reductions $P^\cdot \otimes \Lambda _n$. For a $\Lambda $--module $M$ and an integer $k$ let $M[\lambda ^k]$ denote the submodule of $M$ of elements annihilated by $\lambda ^k$.  Then from the exact sequence
$$
\begin{CD}
0@>>> P^\cdot @>\lambda ^k>> P^\cdot @>>> P^\cdot /\lambda ^k@>>> 0
\end{CD}
$$
one obtains for every $n$ a short exact sequence
$$
0\rightarrow H^i(P^\cdot)\otimes \Lambda _n\rightarrow H^i(K_n)\rightarrow H^{i+1}(P^\cdot)[\lambda ^n]\rightarrow 0.
$$
These short exact sequences fit together to form a short exact sequence of projective systems, where the transition maps
$$
H^{i+1}(P^\cdot)[\lambda ^{n+1}]\rightarrow H^{i+1}(P^\cdot)[\lambda ^n]
$$
are given by multiplication by $\lambda $.  Since $H^{i+1}(P^\cdot)$ is of finite type and in particular has bounded $\lambda $--torsion, it follows that the map of projective systems
$$
H^i(P^\cdot)\otimes \Lambda _n\rightarrow H^i(K_n)
$$
has AR-null kernel and cokernel. This proves (i).

For (ii), note that if $z\in P^i$ is a closed element then modulo $\lambda $ the element $z$ is a boundary.  Write $z = \lambda z'+d(a)$ for some $z'\in P^i$ and $a\in P^{i-1}$.  Since $P^{i+1}$ is flat over $\Lambda $ the element $z'$ is closed.  It follows that $H^i(P^\cdot) = \lambda H^i(P^\cdot)$.  Since $H^i(P^\cdot)$ is a finitely generated $\Lambda $--module, Nakayama's lemma implies that $H^i(P^\cdot) = 0$.  Thus by (i) the system $H^i(K_n)$ is AR--isomorphic to $0$ which implies (ii).
\end{proof}

\begin{thm} The functor $F$ in \ref{comparisonfunctor} is an equivalence of triangulated categories.
\end{thm}
\begin{proof}
Since the $\Ext^{-1}$'s involved in \ref{cor-compa-deligne} are finite dimensional for bounded constructible complexes, the full faithfulness follows from \ref{cor-compa-deligne}.

For the essential surjectivity, note first that any object $K\in \D^b_{c, \text{Del}}(X, \Lambda )$ is induced by a complex $M\in \D_c(X^\NN, \Lambda _\bullet )$ by restriction.  For example represent each $K_n$ by a homotopically injective complex $I_n$ in which case the morphisms $K_{n+1}\rightarrow K_n$ defined in the derived category can be represented by actual maps of complexes $I_{n+1}\rightarrow I_n$.  By \ref{normalchar} the complex $M$ is normalized and by the preceding lemma the corresponding object of $\DD_c(X, \Lambda )$ lies in $\DD_c^b(X, \Lambda )$.  It follows that if $\overline M\in \DD_c^b(X, \Lambda )$ denotes the image of $M$ then $K$ is isomorphic to $F(\overline M)$.
\end{proof}

\begin{rem} One can also define categories $\DD_c(\X, \mathbb{Q}_l)$.  There are several different possible generalizations of the classical definition of this category for bounded complexes on noetherian schemes.  The most useful generalizations seems to be to consider the full subcategory $\mc T$ of $\DD_c(\X, \Z_l)$ consisting of complexes $K$ such that for every $i$ there exists an integer $n\geq 1$ such that $\mc H^i(K)$ is annihilated by $l^n$.   Note that if $K$ is an unbounded complex there may not exist an integer $n$ such that $l^n$ annihilates all $\mc H^i(K)$.  Furthermore, when $\X$ is not quasi--compact the condition is \emph{not} local on $\mc X$. Nonetheless, by \cite[2.1]{Nee01} we can form the quotient of $\DD_c(\X, \mathbb{Z}_l)$ by the subcategory $\mc T$ and we denote the resulting triangulated category with $t$--structure (induced by the one on $\DD_c(\X, \mathbb{Z}_l)$)  by $\DD_c(\X, \mathbb{Q}_l)$. If $\X$ is quasi--compact and $F, G\in \DD^b(\X, \mathbb{Z}_l)$ one has
$$
\text{Hom}_{\DD^b(\X, \mathbb{Z}_l)}(F, G)\otimes \Q\simeq \text{Hom}_{\DD^b(\X, \mathbb{Q}_l)}(F, G).
$$

Using a similar $2$--categorical limit method as in \cite[1.1.3]{WeilII} one can also define a triangulated category $\DD_c(\X, \overline {\mathbb{Q}}_l).$
\end{rem}

%\begin{rem} For example, if $S$ is the spectrum of a finite
%field  and $\X$ is a $k$-variety,  then the $\Ext^{-1}$'s involved
%are of finite dimension for bounded constructible complexes. It
%follows that the natural functor (cf.~\ref{hatMbounded})
%$M\ra(\hat M)_n$ is a fully faithful functor from $\DD^b_c(\X)$ to
%Deligne's  bounded derived category of constructible $l$-adic
%sheaves \cite{WeilII}.
%end{rem}

\section{$\Rhom$} We define the bifunctor
$$\Rhom:\DD_c(\X)^{\mathrm{opp}}\times\DD_c(\X)\ra\D(\X)$$ by the formula
$$\Rghom(M,N)=\Rhom_{\Lambda_\bullet}(\hat M,\hat N).$$

Recall that $\D_c(\X,\Lambda_n)$ denotes the usual derived
category of complexes of $\Lambda_n$-modules with constructible
cohomology.
\begin{prop}\label{lemme-normal-chapeau} Let $M\in\DD_c^-(\X)$ and $N\in\DD_c^+(\X)$, then
    $\Rghom(M,N)$ has constructible cohomology and is normalized. Therefore, it defines an additive
functor
$$\Rghom:\DD^-_c(\X)^{\textup{opp}}\times\DD_c^+(\X)\ra\DD_c^+(\X).$$
\end{prop}
\begin{proof} One can assume $M,N$ normalized. By~\ref{calc-loc}, one has the formula
\begin{equation}\label{locrhom}e_n^{-1}\Rhom(M,N)=\Rhom(e_n^{-1}M,e_n^{-1}N).\end{equation}
>From this it follows that $\Rghom(M, N)$ has constructible cohomology.

 By~\ref{locrhom} and~\ref{normalchar}, to prove that $\Rghom(M, N)$ is normalized we have to show that
$$\Lambda_n\Otimes_{\Lambda_{n+1}}\Rhom_{\Lambda_{n+1}}
(\Lambda_{n+1}\Otimes_\Lambda{R}\pi_*M,\Lambda_{n+1}\Otimes_\Lambda{R}\pi_*N)\ra
\Rhom_{\Lambda_{n}}
(\Lambda_{n}\Otimes_\Lambda{R}\pi_*M,\Lambda_{n}\Otimes_\Lambda{R}\pi_*N)$$
is an isomorphism.
By~\ref{calcul-pou-rhom}, both
$M_{n+1}=\Lambda_{n+1}\Otimes_\Lambda{R}\pi_*M$ and
$N_{n+1}=\Lambda_{n+1}\Otimes_\Lambda{R}\pi_*N$ are constructible
complexees of $\Lambda_{n+1}$ sheaves on $U_\et$. One is reduced to the formula
$$\Lambda_n\Otimes_{\Lambda_{n+1}}\Rhom_{\Lambda_{n+1}}
(M_{n+1},N_{n+1})\ra \Rhom_{\Lambda_{n}}
(\Lambda_{n}\Otimes_{\Lambda_{n+1}}M_{n+1},\Lambda_{n}\Otimes_{\Lambda_{n+1}}N_{n+1})$$
for  our constructible complexes $M,N$ on $U_\et$. This assertion is
well-known (and is easy to prove), cf.~\cite[lemma II.7.1, II.7.2]{SGA5}. \end{proof}

\begin{rem}\label{rem-rhom} Using almost the same proof, one can define a functor
$$\Rghom:\DD^b_c(\X)^{\mathrm{opp}}\times\DD_c(\X)\ra\DD_c(\X).$$
\end{rem}

\section{$\RgHom$}
Let $M,N$ in $\DD_c^-(\X),\DD^+_c(\X)$ respectively.  We define
the functor
$$\RgHom:\DD_c^{-\textup{opp}}(\X)\times\DD^+_c(\X)\ra\textup{Ab}$$ by the formula

\begin{equation}\label{RHom}
    \RgHom(M,N)=\RHom_{\Lambda_\bullet}(\hat M,\hat
N).
\end{equation}
 By~\ref{small-hom}, one has
 $$H^0\RgHom(M,N)=\Hom_{\D_c(\X^\NN)}(\hat M,\hat N)=\Hom_{\DD_c(\X)}(M,N).$$
One has
$$\RHom_{\Lambda_\bullet}(\hat M,\hat
N)=\RHom_{\Lambda_\bullet}(\Lambda_\bullet,\Rhom(\hat M,\hat
N)).$$ By~\ref{lemme-normal-chapeau}, $\Rhom(\hat M,\hat N)$ is
constructible and normalized. Taking $H^0$, we get the formula
$$\Hom_{\D_c(\X^\NN)}(\hat M,\hat
N)=\Hom_{\D_c(\X^\NN)}(\Lambda_\bullet,\Rhom(\hat M,\hat N)).$$
By~\ref{small-hom}, we get therefore
the formula\begin{equation}\label{RHom=Hom}
    \Hom_{\DD_c(\X)}(M,N)=\Hom_{\DD_c(\X)}(\Lambda,\Rghom(M,N))
\end{equation}

In summary, we have gotten the following result.

\begin{prop}\label{Rhom-versus-RHom}Let $M,N$ in $\DD_c^-(\X),\DD^+_c(\X)$
respectively.One has
$$\Hom_{\DD_c(\X)}(M,N)=H^0\RgHom(M,N)=\Hom_{\DD_c(\X)}(\Lambda,\Rghom(M,N)).$$
\end{prop}

\begin{rem}\label{def-ext} Accordingly, one defines $$\egxt^{*}(M,N)=\mc
H^*(\Rghom(M,N))\text{ and }\Egxt^*(M,N)=\mc H^*(\RgHom(M,N))$$
and
$$\Hgom(M,N)=Hom_{\DD_c(\X)}(M,N)=\mc H^0(\RgHom(M,N)).$$\end{rem}

%\begin{rem}[For Martin] Under the assumptions of~\ref{Rhom-versus-RHom},
%$\Rghom(M,N)\in\DD^+$ if $\X$ is of finite type over $S$. But it's
%ot clear a priori for $\X$ only locally of finite type.\end{rem}

\section{Tensor product}

Let $M,N\in\DD_c(\X)$. We define the total tensor product
$$M\Otimes_{\Lambda}N=\hat M\Otimes_{\Lambda_\bullet}\hat N.$$
It defines a bifunctor
$$\DD_c(\X)\times\DD_c(\X)\ra\D_{\textrm{cart}}(\X).$$

\begin{prop}\label{homandtensor} For any $L, N, M \in \DD_c(\X, \Lambda )$ we have
\begin{equation*}
\Rghom (L\Otimes N, M)\simeq \Rghom (L, \Rghom (N, M)).
\end{equation*}
\end{prop}
\begin{proof}
By definition this amounts to the usual adjunction formula
\begin{equation*}
\Rhom (\hat L\Otimes \hat N, \hat M)\simeq \Rhom (\hat L, \Rhom (\hat N, \hat M)).
\end{equation*}
\end{proof}

\begin{cor} For any $L, M\in \DD_c(\X, \Lambda )$ there is a canonical evaluation morphism
\begin{equation*}
\text{\rm ev}:\Rghom (L, M)\Otimes L\rightarrow M.
\end{equation*}
\end{cor}
\begin{proof}
The morphism $\text{\rm ev}$ is defined to be the image of the identity map under the isomorphism
\begin{equation*}
\Rghom (\Rghom (L, M), \Rghom (L, M))\simeq \Rghom (\Rghom (L, M)\Otimes L, M)
\end{equation*}
provided by \ref{homandtensor}.
\end{proof}

\section{Duality} Let's denote by $f:\X\ra S$ the structural
morphism. Let
\begin{equation}\label{omegan}
   \Omega_n=f^!\Lambda_n(\dim(S))[2\dim(S)]
\end{equation}
be the (relative) dualizing complex of $\X$ (ringed by
$\Lambda_n$). Notice that $f^*\Lambda_n=\Lambda_n$, with a slight
abuse of notation.

\subsection{Construction of the dualizing complex}
\begin{prop}\label{Omega-adic} One has
$\Omega_n=\Lambda_n\Otimes_{\Lambda_{n+1}}\Omega_{n+1}$.
\end{prop}
\begin{proof}
The key point is the following lemma:
\begin{lem} Let $M$ be a complex of sheaves of $\Lambda _{n+1}$--modules of finite injective dimension.  Then there is a canonical isomorphism
\begin{equation*}
M\Otimes _{\Lambda _{n+1}}\Lambda _n\simeq \Rhom _{\Lambda _{n+1}}(\Lambda _n, M).
\end{equation*}
\end{lem}
\begin{proof}
Let $S$ denote the acyclic complex on $\X$
$$\cdots\Lambda_{n+1}\stackrel{l^{n+1}}{\ra}\Lambda_{n+1}
\stackrel{l}{\ra}\Lambda_{n+1}\stackrel{l^{n+1}}{\ra}\Lambda_{n+1}\ra
\cdots,
$$
where the map $S^{i}\rightarrow S^{i+1}$ is given by multipliciation by $l$ if $i$ is even and multiplication by $l^{n+1}$ if $i$ is odd.
Let $P$ denote the truncation $\sigma _{\leq 0}S$ (the terms of $S$ in degrees $\leq 0$). Then $P$ is a projective resolution of $\Lambda _n$ viewed as a $\Lambda _{n+1}$--module and $\widehat P:= \text{Hom}^\cdot (P, \Lambda _{n+1})$ is isomorphic to $\sigma _{\geq 1}S$ and is also a resolution of $\Lambda _n$.  The diagram
\begin{equation*}
\begin{CD}
\cdots \rightarrow \Lambda _{n+1}@>l^{n+1}>> \Lambda _{n+1}@>>> 0 \\
@VVV  @V\times lVV @VVV \\
0@>>> \Lambda _{n+1}@>l^{n+1}>> \Lambda _{n+1} \rightarrow \cdots
\end{CD}
\end{equation*}
defines a morphism of complexes $P\rightarrow \widehat P[1]$ whose cone is quasi--isomorphic to $S[1]$. Let $M$ be a boundex complex of injectives.  We then obtain a morphism
\begin{equation}\label{hommap}
M\Otimes _{\Lambda _{n+1}}\Lambda _n\simeq M\otimes P\simeq \text{Hom}^\cdot (\widehat P, M)\rightarrow \text{Hom}^{\cdot }(P, M)\simeq \Rhom _{\Lambda _{n+1}}(\Lambda _n, M).
\end{equation}
The cone of this morphism is isomorphic to $\Rhom (S, M)$ which is zero since $S$ is acyclic. It follows that \ref{hommap} is an isomorphism.
\end{proof}

In particular
$$\Rhom_{\Lambda_{n+1}}(\Lambda_n,f^!\Lambda_{n+1})=\Lambda_n\Otimes_{\Lambda_{n+1}}f^{!}(\Lambda_{n+1}).$$
On the other hand, by \cite[4.4.3]{Las-Ols05}, one has
$$\Rhom_{\Lambda_{n+1}}(\Lambda_n,f^!\Lambda_{n+1})=
f^!\Rhom_{\Lambda_{n+1}}(\Lambda_n,\Lambda_{n+1})=f^!\Lambda_n.$$
Twisting and shifting, one gets \ref{Omega-adic}.\end{proof}

\begin{pg}\label{LETdata} Let $U\rightarrow \X$ be an object of $\LE(\X)$ and let $\epsilon :\X|_U\rightarrow U_\et$ be the natural morphism of topos. Let us describe more explicitly the morphism $\epsilon$. Let $\LE(\mc X)_{|U}$ denote the category of morphisms $V\rightarrow U$ in $\LE (\mc X)$.  The category $\LE(\mc X)_{|U}$ has  a Grothendieck topology induced by the topology on $\LE(\mc X)$, and the resulting topos is
canonically isomorphic to the localized topos $\X_{\le|U}$.  Note
that there is a natural inclusion $\LE(U)\hookrightarrow \LE(\mc
X)_{|U}$ but this is not an equivalence of categories since for an
object $(V\rightarrow U)\in \LE(\mc X)_{|U}$ the morphism
$V\rightarrow U$ need not be smooth. Viewing $\X_{\le|U}$ in this
way, the functor $\epsilon^{-1}$ maps $F$ on $U_\et$ to
$F_V=\pi^{-1}F\in V_\et$ where $\pi:V\ra U\in\LE(\X)_{|U}$. For a
sheaf $F\in \X_{\le |U}$ corresponding to a collection of sheaves
$F_V$, the sheaf $\epsilon _*F$ is simply the sheaf $F_U$.

In particular, the functor $\epsilon _*$ is exact and,
accordingly $H^*(U,F)=H^*(U_{\et},F_U)$ for any sheaf of
$\Lambda$ modules of $\X$.
\end{pg}

\begin{thm}\label{def-Omega-adic} There exists a  normalized complex $\Omega_\bullet\in\D_c(\X^\NN)$,
unique up to canonical isomorphism, inducing the $\Omega _n$.\end{thm} \begin{proof} The
topos $\X_\le^\NN$ can be described by the site $\mc S$ whose
objects are pairs $(n,u:U\ra\X)$ where $u$   is a lisse-{\'e}tale
open and $n\in\NN$. We want to use the gluing theorem~\cite[2.3.3]{Las-Ols05}.

\begin{itemize}
    \item Let us describe the localization morphisms
    explicitly. Let $(U,n)$ be in $\mc S$. An object of the localized topos $\X^\NN_{|(U,n)}$
    is equivalent to giving for every $U$--scheme of
finite type $V\rightarrow U$, such that the composite
$\alpha:V\rightarrow U\rightarrow \X$ is smooth of relative
dimension $d_\alpha$, a projective system
$$F_V=(F_{V,m},m\leq n)$$
where $F_{V,m}\in V_{\et }$ together with morphisms
$f^{-1}F_V\rightarrow F_{V'}$ for $U$--morphisms $f:V'\rightarrow
V$. The localization morphism
$${j_n}:\X^\NN_{|(U,n)}\ra \X^\NN$$
is defined by the truncation
$$({j_n}^{-1}F_\bullet)_V=(F_{m,V})_{m\leq n}.$$
We still denote ${j_n}^{-1}\Lambda_\bullet=(\Lambda_m)_{m\leq n}$
by $\Lambda_\bullet$ and we ring $\X^\NN_{|(U,n)}$ by
$\Lambda_\bullet$.

\item Notice that $\pi:\X^\NN\ra\X$ induces
$$\pi_n:\X^\NN_{|(U,n)}\ra\X_{|U}$$ defined by $\pi_n^{-1}(F)=(F)_{m\leq n}$ (the
constant projective system). One has $$\pi_{n*}(F_m)_{m\leq
n}=\varprojlim_{m\leq n} F_m=F_n.$$

\item As in the proof of \ref{calc-loc}, the morphisms $e_i:\X\ra\X^\NN, i\leq n$ can be localized in
    ${\tilde e}_i:\X_{|U}\ra\X^\NN_{|(U,n)}$, characterized by ${\tilde e}_i^{-1}(F_m)_{m\leq
    n}=F_i$.They form a conservative sets of functors.
\item
    One has a commutative diagram of topos
\begin{equation}\label{diag-comm}
    \xymatrix{
    \X_{|U}\ar[r]^{{\tilde e}_n}\ar[rdd]_\epsilon&\X^\NN_{|(U,n)}
    \ar[d]^{\pi_n}\ar@/^2pc/[dd]^{p_n}\ar[r]^{j_n}&\X^\NN\\
&\X_{|U}\ar[d]^\epsilon\\&U_\et}
\end{equation} One has
$\pi_n^{-1}(\Lambda_n)=(\Lambda_n)_{m\leq n}$ -the constant
projective system with value $\Lambda_n$- which maps to
$(\Lambda_m)_{m\leq n}$ : we will ring $\X_{|U}$ (and also both
$\X$ and $U_\et$) by $\Lambda_n$ and therefore the previous
diagram is a diagram of ringed topos. Notice that $e_n^{-1}=e_n^*$
implying the exactness of $e_n^*$.
 \item Let us define \begin{equation}\label{deflocomega}
 \Omega_{U,n}=L\pi_n^*\Omega_{n|U}=Lp_n^*K_{U,n}\langle-d_\alpha\rangle\end{equation}
    where $K_{U,n}\in\D_c(U_\et,\Lambda_n)$ is the dualizing
    complex.

    \item Let $f:(V,m)\ra (U,n)$ be a morphism in $\mc S$. It
    induces a commutative diagram of ringed topos
    $$\xymatrix{\X^\NN_{|(V,m)}\ar[d]^{p_m}\ar[r]^f&\X^\NN_{|(U,n)}\ar[d]^{p_n}\\
    V_\et\ar[r]^f&U_\et}.$$ By the construction of the dualizing complex in~\cite{Las-Ols05}
    and~\ref{Omega-adic}, one has therefore
   \begin{equation}\label{Klocal}
   Lf^*\Omega_{U,n}=L\pi_{m}^*(\Lambda_m\Otimes_{\Lambda_n}\Omega_{n|V})=L\pi_m^*\Omega_{m|V}=
   \Omega_{V,m}.
   \end{equation}
   Therefore, $\Omega_{U,n}$ defines locally an object $\D_c(\mc
S,\Lambda_\bullet)$. Let's turn to the
$\ext$'s.
  \item The morphism of topos $\pi_n:\X^\NN_{|(U,n)}\ra \X_{|U}$
    is defined by $\pi_n^{-1}F=(F)_{m\leq n}$. One
    has therefore $$\pi_{n*}F=F_n\text{ and } p_{n*}F=F_{n,U}.$$
    In particular, one gets the exactness of $p_{n*}$
    and the formulas  \begin{equation}\label{pi=en}Rp_{n*}=p_{n*}\text{ and
    }\pi_{n*}={\tilde e}_n^*.\end{equation}
Using (\ref{diag-comm}) we get the formula
    \begin{equation}\label{p=en}p_{n*}Lp_n^*=\epsilon_*{\tilde e}_n^*Lp_n^*=\epsilon_*\epsilon^*=\Id.\end{equation}

 Therefore
 one has
\begin{eqnarray*}
\Ext^i(Lp_n^*K_{U,n},Lp_n^*K_{U,n})& = & \Ext^i(K_{U,n},p_{n*}Lp_n^*K_{U,n})\\
&=& \Ext^i(K_{U,n},K_{U,n}) \text{ \ \ \ by \ref{p=en} }\\
& = & H^i(U_\et,\Lambda_n) \text{ \ \ \ by duality}.
\end{eqnarray*}
 By sheafification, one gets
 $$\ext^i(Lp_n^*K_{U,n},Lp_n^*K_{U,n})=\Lambda_\bullet\text{ for
 }i\not=0\text{ and }\Lambda_\bullet\text{ else.}$$

Therefore, the local data $(\Omega_{U,n})$ has vanishing negative
$\ext$'s. By \cite[3.2.2]{Las-Ols05}, there exists a unique
$\Omega_\bullet\in\D_c(\X,\Lambda)$ inducing $\Omega_{U,n}$ on
each $\X^\NN_{|(U,n)}$ . By~\ref{deflocomega}, one has

\item Using the formula $j_n\circ {
e}_n=\tilde e_n\circ j$ (\ref{diag-comm}) and \ref{deflocomega}, one
obtains
$$(e_n^*\Omega_\bullet)_{|U}=e_n^*\Omega_{(n,U)}=\Omega_{n|U}.$$
By~\cite[3.2.2]{Las-Ols05}, the isomorphisms glue to define a functorial
isomorphism
$$e_n^*\Omega_\bullet=\Omega_{n}.$$ By~\ref{Omega-adic} and \ref{normalchar},
$\Omega_\bullet$ is normalized with constructible cohomology.

\item The uniqueness is a direct consequence
of~(\ref{compa-deligne}).\end{itemize}

\end{proof}

\subsection{The duality theorem} Let $M$ be a normalized complex.
 By \ref{calc-loc}, one has
\begin{equation}\label{loc-D}
    e_n^{-1}\Rhom(M,\Omega)\ra\Rhom(e_n^{-1}M,e_n^{-1}\Omega)
\end{equation}
(\ref{calc-loc}). The complex $\Omega $ is of locally finite quasi--injective dimension in the following sense.  If $\X$ is quasi--compact, then each
 $\Omega_n$ is of finite quasi-injective
dimension, bounded by some integer $N$ depending only on
$\X$ and $\Lambda$, but not $n$.  Therefore in the quasi--compact case one has $$\egxt^i(M,\Omega)=0\text{ for any
}M\in\DD^{\geq 0}_c(\X)\text{ and }i\geq N.$$

Let's now prove the duality theorem.
\begin{thm}\label{dual-adic} Let
$D:\DD_c(\X)^{\textup{opp}}\ra\D(\X)$ be the functor defined by
$D(M)=\Rghom(M,\Omega)=\Rhom_{\Lambda_\bullet}(\hat
M,\hat\Omega)$.
\begin{enumerate}
    \item The essential image of $D$ lies in $\D_c(\A)$.
    \item If $D:\DD_c(\X)^{\textup{opp}}\ra\DD_c(\X)$  denotes the induced
functor, then $D$ is involutive and maps $\DD_c^-(\X)$ into $\DD^+_c(\X)$.
\end{enumerate}
\end{thm}

\begin{proof} Both assertions are local on $\X$ so we may assume that $\X$ is quasi--compact.  Because $\Omega $ is of finite quasi--injective dimension, to prove the first point it suffices to prove (1) for bounded below complexes.  In this case the result follows from \ref{lemme-normal-chapeau}.

For the second point, one can assume $M$ normalized (because
$\hat M$ is constructible (\ref{th-normalisation}) and normalized.
Because $\Omega$ is normalized (\ref{Omega-adic}), the
tautological biduality morphism
$$\Rhom(M,\Rhom(M,\Omega),\Omega)\ra M$$ defines a morphism
$$DD(M)\ra M.$$ Using~\ref{loc-D}, one is reduced to the analogous formula
$$D_nD_n(e_n^{-1}M)=e_n^{-1}M$$ where $D_n$ is the dualizing functor on $\D_c(\X_\le,\Lambda_n)$,
which is proven in \textit{loc. cit.}.\end{proof}

\begin{cor}\label{7.8} For any $N, M\in \DD_c(\X, \Lambda )$ there is a canonical isomorphism
\begin{equation*}
\Rghom (M, N)\simeq \Rghom (D(N), D(M)).
\end{equation*}
\end{cor}
\begin{proof}
Indeed by \ref{homandtensor} we have
\begin{equation*}
\Rghom (D(N), D(M))\simeq \Rghom (D(N)\Otimes M, \Omega )\simeq \Rghom (M, DD(N))\simeq \Rghom (M, N).
\end{equation*}
\end{proof}

\section{The functors $Rf_*$ and $Lf^*$}

\begin{lem} Let $f:\X\rightarrow \mc Y$ be a morphism of finite type between  $S$--stacks. Then for any integer $n$ and $M\in \D_c^+(\X, \Lambda _{n+1})$ the natural map
\begin{equation*}
Rf_*M\Otimes _{\Lambda _{n+1}}\Lambda _n\rightarrow Rf_*(M\Otimes _{\Lambda _{n+1}}\Lambda _n)
\end{equation*}
is an isomorphism.
\end{lem}
\begin{proof}
The assertion is clearly local in the smooth topology on $\mc Y$ so we may assume that $\mc Y$ is a scheme.  Furthermore, if $X_\bullet \rightarrow \X$ is a smooth hypercover by schemes and $M_\bullet \in \D_c(X_\bullet, \Lambda _{n+1})$ is the complex corresponding to $M$ under the equivalence of categories $\D_c(X_{\bullet, \et }, \Lambda _{n+1})\simeq \D_c(\X, \Lambda _{n+1})$ then by \cite[9.8]{Ols05} it suffices to show the analogous statement for the morphism of topos
$$
f_{\bullet }:X_{\bullet, \et }\rightarrow \mc Y_{\text{et}}.
$$
Furthermore by a standard spectral sequence argument (using the sequence defined in \cite[9.8]{Ols05}) it suffices to prove the analogous result for each of the morphisms $f_n:X_{n, \et }\rightarrow \mc Y_{\et }$, and hence it suffices to prove the lemma for a finite type morphism of schemes of finite type over $S$ with the \'etale topology where it is standard.
\end{proof}

\begin{prop} Let $M = (M_n)_{n}$ be a bounded below $\lambda $--complex on $\X$.  Then for any integer $i$ the system $R^if_*M= (R^if_*M_n)_n$ is almost adic.
\end{prop}
\begin{proof}
The assertion is clearly local on $\mc Y$, and hence we may assume that both $\mc X$ and $\mc Y$ are quasi--compact.

By the same argument proving \cite[9.10]{Ols05} and \cite[Th. finitude]{SGA4.5}, the sheaves $R^if_*M_n$ are constructible. The result then follows from \cite[V.5.3.1]{SGA5} applied to the category of constructible sheaves on $\X_{\le }$.
\end{proof}

Now consider the morphism of topos $f_\bullet :\X^{\mathbb{N}}\rightarrow \mc Y^\N$ induced by the morphism $f$.  By the above, if $M\in \D^+(\X^{\mathbb{N}})$ is a $\lambda $--complex then $Rf_*M$ is a $\lambda $--complex on $\mc Y$.  We therefore obtain a functor
\begin{equation*}
Rf_*:\DD_c^+(\X, \Lambda )\rightarrow \DD_c^+(\mc Y, \Lambda ).
\end{equation*}

It follows immediately from the definitions that the pullback functor $Lf^*:\D_c(\mc Y^{\mathbb{N}}, \Lambda )\rightarrow \D_c(\mc X^\N, \Lambda )$ take $\lambda $--complexes to $\lambda $--complexes and AR--null complexes to AR--null complexes and therefore induces a functor
\begin{equation*}
Lf^*:\DD_c(\mc Y, \Lambda )\rightarrow \DD_c(\mc X, \Lambda ), M\mapsto Rf_*\hat M.
\end{equation*}

\begin{prop}\label{8.3} Let $M\in \DD_c^+(\X, \Lambda )$ and $N\in \DD_c^-(\mc Y, \Lambda )$. Then there is a canonical isomorphism
\begin{equation*}
Rf_*\Rghom (Lf^*N, M)\simeq \Rghom (N, Rf_*M).
\end{equation*}
\end{prop}
\begin{proof}
We can rewrite the formula as
\begin{equation*}
Rf_*\Rhom (Lf^*\hat N, \hat M)\simeq \Rhom (\hat N, Rf_*\hat M)
\end{equation*}
which follows from the usual adjunction between $Rf_*$ and $Lf^*$.
\end{proof}

\section{The functors $Rf_!$ and $Rf^!$}

\subsection{Definitions}

Let $f:\X\rightarrow \mc Y$ be a finite type morphism of $S$--stacks and let $\Omega _{\mc X}$ (resp. $\Omega _{\mc y}$) denote the dualizing complex of $\X$ (resp. $\mc Y$).  Let $D_{\X}:\DD_c(\X)\rightarrow \DD_c(\X)$ denote the functor $\Rghom (-, \Omega _{\X}):\DD_c(\X, \Lambda )\rightarrow \DD_c(\X, \Lambda )$ and let $D_{\mc Y}:\DD_c(\mc Y, \Lambda )\rightarrow \DD_c(\mc Y, \Lambda )$ denote $\Rghom (-, \Omega _{\mc Y})$.
We then define
\begin{equation*}
Rf_!:= D_{\mc Y}\circ Rf_*\circ D_{\X}:\DD_c^-(\X, \Lambda )\rightarrow \DD_c^-(\mc Y, \Lambda )
\end{equation*}
and
\begin{equation*}
Rf^!:= D_{\mc X}\circ Lf^*\circ D_{\mc Y}:\DD_c(\mc Y, \Lambda )\rightarrow \DD_c(\X, \Lambda ).
\end{equation*}

\begin{lem} For any $N\in \DD_c^-(\mc X, \Lambda )$ and $M\in \DD_c^+(\mc Y, \Lambda )$ there is a canonical isomorphism
\begin{equation*}
Rf_*\Rghom (N, Lf^!M)\simeq \Rghom (Rf_!N, M).
\end{equation*}
\end{lem}
\begin{proof}
Set $N' = D_\X(N)$ and $M':= D_{\mc Y}(M)$.  Then by \ref{7.8} the formula can be written as
\begin{equation*}
Rf_*\Rghom (Lf^*M', N')\simeq \Rghom (M', Rf_*N')
\end{equation*}
which is \ref{8.3}.
\end{proof}

\begin{lem}\label{9.3} If $f$ is a smooth morphism of relative dimension $d$, then there is a canonical isomorphism $Rf^!(F)\simeq f^*F\langle d\rangle $.
\end{lem}
\begin{proof} By the construction of the dualizing complex and \cite[4.5.2]{Las-Ols05}  we have $\Omega _{\mls X} \simeq f^*\Omega _{\mls Y}\langle d\rangle $.  From this and biduality \ref{dual-adic} the lemma follows.
\end{proof}

If $f$ is a closed immersion, then we can also define the functor of sections with support $\underline {H}_\X^0$ on the category of $\Lambda _\bullet $--modules in $\mls Y^{\mathbb{N}}$. This functor is right adjoint to $f_*$ and taking derived functors we obtain an adjoint pair of functors
$$
Rf_*:\D_c(\mls X^{\mathbb{N}}, \Lambda _\bullet )\rightarrow \D_c(\mls Y^{\mathbb{N}}, \Lambda _\bullet )
$$
and
$$
R\underline {H}_\X^0:\D_c(\mls Y^{\mathbb{N}}, \Lambda _\bullet )\rightarrow \D_c(\mls X^{\mathbb{N}}, \Lambda _\bullet ).
$$
Both of these functors take $AR$--null complexes to $AR$--nul complexes and hence induce adjoint functors on the categories $\DD_c(\mls Y, \Lambda )$ and $\DD_c(\mls X, \Lambda )$.

\begin{lem} If $f$ is a closed immersion, then $\Omega _{\mls X} = f^*R\underline {H}_\X^0\Omega _{\mls Y}$.
\end{lem}
\begin{proof} By the gluing lemma this is a local assertion in the topos $\mls X^{\mathbb{N}}$ and hence the result follows from \cite[4.6.1]{Las-Ols05}.
\end{proof}

\begin{prop} If $f$ is a closed immersion, then $f^! = f^*R\underline {H}_\X^0$ and $Rf_* = Rf_!$.
\end{prop}
\begin{proof} This follows from the same argument proving \cite[4.6.2]{Las-Ols05}.
\end{proof}

Finally using the argument of \cite[4.7]{Las-Ols05} one shows:
\begin{prop}\label{9.6} If $f$ is a universal homeomorphism then $f^*\Omega _{\mls X} = \Omega _{\mls Y}$, $Rf^! = f^*$, and $Rf_! = Rf_*$.
\end{prop}

There is also a projection formula
\begin{equation}\label{projection}
Rf_!(A\Otimes f^*B)\simeq Rf_!A\Otimes B
\end{equation}
for $B\in \DD_c^-(\mls Y, \Lambda )$ and $A\in \DD_c(\mls X, \Lambda )$.
This is shown by  the same argument used to prove \cite[4.4.2]{Las-Ols05}.

\section{Computing cohomology using hypercovers}

For this we first need some cohomological descent results.

Let $\mls X$ be a algebraic stack over $S$ and $X_\bullet
\rightarrow \mls X$ a strictly simplicial smooth hypercover with
the $X_i$ also $S$--stacks.  We can then also consider the
topos of projective systems in $X_{\bullet, \le}$ which we denote
by $X_{\bullet }^{\mathbb{N}}$.

\begin{defn} (i) A sheaf $F$ of $\Lambda _\bullet $--modules in $X_{\bullet }^{\mathbb{N}}$ is \emph{almost adic} if it is cartesian and if for every $n$ the restriction $F|_{X_{n, \le }}$ is almost adic.

(ii) An object $C\in \D(X_\bullet ^{\mathbb{N}}, \Lambda _\bullet )$ is a \emph{$\lambda $--complex} if for all $i$ the cohomology sheaf $\mc H^i(C)$ is almost adic.

(iii) An object $C\in \D(X_\bullet ^{\mathbb{N}}, \Lambda _\bullet )$ is \emph{almost zero} if for every $n$ the restriction of $C$ to $X_n$ is almost zero.

(iv) Let $\D_c(X_{\bullet }^{\mathbb{N}}, \Lambda _\bullet )\subset \D(X_{\bullet }^{\mathbb{N}}, \Lambda _\bullet )$ denoted the triangulated subcategory whose objects are the $\lambda $--complexes. The category $\DD_c(X_{\bullet }, \Lambda )$ is the quotient of $\D_c(X_{\bullet }^{\mathbb{N}}, \Lambda _\bullet )$ by the full subcategory of almost zero complexes.
\end{defn}

As in \ref{gen-proj} we have the projection morphism
\begin{equation*}
\pi :(X_{\bullet }^{\mathbb{N}}, \Lambda _\bullet )\rightarrow (X_{\bullet, \le }, \Lambda )
\end{equation*}
restricting for every $n$ to the morphism $(X_{n }^{\mathbb{N}}, \Lambda _\bullet )\rightarrow (X_{n, \le }, \Lambda )$ discussed in \ref{gen-proj}.  By \ref{E1} the functor $R\pi _*:\D(X_{\bullet }^{\mathbb{N}}, \Lambda _\bullet )\rightarrow \D(X_\bullet, \Lambda )$ takes almost zero complexes to $0$.  By the universal property of the quotient category it follows that there is an induced functor
\begin{equation*}
R\pi _*:\DD_c(X_{\bullet } , \Lambda )\rightarrow \D(X_{\bullet }, \Lambda ).
\end{equation*}
We also define a normalization functor
\begin{equation*}
\DD_c(X_{\bullet }, \Lambda )\rightarrow \D(X_{\bullet }^{\mathbb{N}}, \Lambda _{\bullet }), M\mapsto \hat M
\end{equation*}
by setting $\hat M:= L\pi ^*R\pi _*(M)$.

\begin{prop}\label{descentprop} Let $M\in \D_c(X_{\bullet }^{\mathbb{N}}, \Lambda _\bullet )$ be a $\lambda $--complex.  Then $\hat M$ is in $\D_c(X_{\bullet }^{\mathbb{N}}, \Lambda _\bullet )$ and the canonical map $\hat M\rightarrow M$ has almost zero cone.
\end{prop}
\begin{proof} For any integer $n$, there is a canonical commutative diagram of ringed topos
\begin{equation*}
\begin{CD}
X_{\bullet }^{\mathbb{N}}@>r_n>> X_n^{\bullet }\\
@V\pi VV @VV\pi V \\
X_{\bullet , \le }@>r_n>> X_{n, \le},
\end{CD}
\end{equation*}
where $r_n$ denotes the restriction morphisms.  Furthermore, the functors $r_{n*}$ are exact and take injectives to injectives.  It follows that for any $M\in \D(X_{\bullet }^{\mathbb{N}}, \Lambda _\bullet )$ there is a canonical isomorphism
$$
R\pi _*(r_{n*}(M))\simeq r_{n*}R\pi _*(M).
$$
>From the definition of $\pi ^*$ it also follows that $r_{n*}L\pi ^* = L\pi ^*r_{n*}$, and from this it follows that the restriction of $\hat M$ to $X_n$ is simply the normalization of $M|_{X_n}$.  From this and \ref{th-normalisation} the statement that $\hat M\rightarrow M$ has almost zero cone follows.

To see that $\hat M\in \D_c(X_{\bullet }^{\mathbb{N}}, \Lambda _\bullet )$, note that by \ref{th-normalisation} we know what for any integers $i$ and $n$  the restriction $\mc H^i(\hat M)|_{X_n}$ is a constructible (and in particular cartesian) sheaf on $X_{n, \le}$.  We also know by \ref{commutU} that for any $n$ and smooth morphism $U\rightarrow X_n$, the restriction of $\mc H^i(\hat M)$ to $U_\et $ is equal to $\mc H^i(\widehat {M_U})$.  From this and \ref{th-normalisation} it follows that the sheaves $\mc H^i(\hat M)$ are cartesian.   In fact, this shows that if $\mc F\in \D_c(\X^\NN, \Lambda _\bullet )$ denotes the complex obtain from the equivalence of categories (cohomological descent as in \cite[2.2.3]{Las-Ols05})
$$
\D_c(X_\bullet^\NN, \Lambda _\bullet )\simeq \D_c(\X^\NN, \Lambda _\bullet ),
$$
then $\mc H^i(\hat M)$ is the restriction to $X_\bullet ^\NN$ of the sheaf $\mc H^i(\hat F)$.
\end{proof}

As in \ref{small-hom} it follows that the normalization functor induces a left adjoint to the projection $\D_c(X_{\bullet }^{\mathbb{N}}, \Lambda _\bullet )\rightarrow \DD_c(X_\bullet , \Lambda )$.

Let $\epsilon :X_{\bullet, \le } \rightarrow \mc X_{\le }$ denote the projection, and write also $\epsilon :X_\bullet ^{\mathbb{N}}\rightarrow \mc X^{\mathbb{N}}$ for the morphism on topos of projective systems.  There is a natural commutative diagram of topos
\begin{equation*}
\begin{CD}
X_{\bullet}^{\mathbb{N}}@>\pi >> X_{\bullet, \le }\\
@V\epsilon VV @VV\epsilon V \\
\mc X^{\mathbb{N}}@>\pi  >> \mc X_{\le }.
\end{CD}
\end{equation*}
By \cite[2.2.6]{Las-Ols05}, the functors $R\epsilon _*$ and $\epsilon ^*$ induce an equivalence of categories
\begin{equation*}
\D_c(X_{\bullet }^{\mathbb{N}}, \Lambda _\bullet )\simeq \D_c(\mc X^{\mathbb{N}}, \Lambda _\bullet ),
\end{equation*}
and the subcategories of almost zero complexes coincide under this equivalence.

We therefore have  obtained

\begin{prop}\label{cohdescentadic} Let $\mls X$ be an algebraic stack over $S$ and $X_\bullet
\rightarrow \mls X$ a strictly simplicial smooth hypercover with
the $X_i$ also $S$--stacks. Then, the morphism
\begin{equation*}
R\epsilon_*:\DD_c(X_{\bullet }, \Lambda )\isom \DD_c(\mc X,
\Lambda )
\end{equation*}is an equivalence with inverse $\epsilon^*$.
\end{prop}

Consider next a morphism of nice stacks $f:\X \rightarrow \mc Y$.  Choose a commutative diagram
\begin{equation*}
\begin{CD}
X_\bullet @>\tilde f>> Y_\bullet \\
@V\epsilon _XVV @VV\epsilon _YV \\
\X @>f>> \mc Y,
\end{CD}
\end{equation*}
where $\epsilon _X$ and $\epsilon _Y$ are smooth (strictly simplicial) hypercovers by nice $S$--stacks.  The functors $Rf_*:\D_c(\X^{\mathbb{N}}, \Lambda _\bullet )\rightarrow \D_c(\mc Y^{\mathbb{N}}, \Lambda _\bullet )$ and $R\tilde f_*:\D_c(X_\bullet ^{\mathbb{N}}, \Lambda _\bullet )\rightarrow \D_c(Y_\bullet ^{\mathbb{N}}, \Lambda _\bullet )$ evidently take almost zero complexes to almost zero complexes and therefore induce functors
\begin{equation*}
Rf_*:\DD _c(\X, \Lambda )\rightarrow \DD_c(\mc Y, \Lambda ), \ \ R\tilde f_*:\DD_c(X_\bullet , \Lambda )\rightarrow \DD_c(Y_\bullet , \Lambda ).
\end{equation*}
It follows from the construction that the diagram
\begin{equation*}
\begin{CD}
\DD _c(\X, \Lambda )@>\ref{cohdescentadic}>> \DD_c(X_\bullet , \Lambda )\\
@VRf_*VV @VVR\tilde f_*V \\
\DD _c(\mc Y, \Lambda )@>\ref{cohdescentadic}>> \DD_c(Y_\bullet , \Lambda )
\end{CD}
\end{equation*}
commutes.

\begin{cor} Let $f:\X\rightarrow \Y$ be a morphism of nice $S$--stacks, and let $X_\bullet \rightarrow \X$ be a strictly simplicial smooth hypercover by nice $S$--stacks.  For every $n$, let $f_n:X_n\rightarrow \Y$ be the projection. Then for any $F\in \DD_c^+(\X, \Lambda )$ there is a canonical spectral sequence in the category of $\lambda $--modules
\begin{equation*}
E_1^{pq} = R^qf_{p*}(F|_{X_p})\implies R^{p+q}f_*(F).
\end{equation*}
\end{cor}
\begin{proof} We take $Y_\bullet \rightarrow \mc Y$ to be the constant simplicial topos associated to $\Y$.  Let $F_\bullet $ denote $\epsilon _X^*F$.  We then have
$$
Rf_*(F) = Rf_*R\epsilon _{X*}(F_\bullet ) = R\epsilon _{Y*}R\tilde f_*(F_\bullet ).
$$
The functor $R\epsilon _{Y*}$ is just the total complex functor (which passes to $\DD_c$), and hence we obtain the corollary from the standard spectral sequence associated to a bicomplex.
\end{proof}

\begin{cor} With notation as in the preceding corollary, let $F\in \DD_c^-(\X, \Lambda )$.  Then there is a canonical spectral sequence in the category of $\lambda $--modules
\begin{equation*}
E_1^{pq} = \mls H^q(D_{\Y}(Rf_{p!}(F)))\implies \mls H^{p+q}(D_{\Y}(Rf_!F)).
\end{equation*}
\end{cor}
\begin{proof} Apply the preceding corollary to $D_{\X}(F)$.
\end{proof}

\section{Kunneth formula}

We prove the Kunneth formula using the method of \cite[\S5.7]{Las-Ols05}.

\begin{lem}\label{9.4} For any $P_1, P_2, M_1, M_2\in \DD_c(\X, \Lambda )$ there is a canonical morphism
\begin{equation*}
\Rghom (P_1, M_1)\Otimes \Rghom (P_2, M_2)\rightarrow \Rghom (P_1\Otimes P_2, M_1\Otimes M_2).
\end{equation*}
\end{lem}
\begin{proof}
By \ref{homandtensor} it suffices to exhibit a morphism
\begin{equation*}
\Rghom (P_1, M_1)\Otimes \Rghom (P_2, M_2)\Otimes P_1\Otimes P_2\rightarrow M_1\Otimes M_2
\end{equation*}
which we obtain from the two evaluation morphisms
\begin{equation*}
\Rghom (P_i, M_i)\otimes P_i\rightarrow M_i.
\end{equation*}
\end{proof}

Let $\mls Y_1$ and $\mls Y_2$ be nice stacks, and set $\mls Y:= \mls Y_1\times \mls Y_2$ with projections $p_i:\mls Y\rightarrow \mls Y_i$. For $L_i\in \DD_c(\mls Y_i, \Lambda )$ let $L_1\Otimes _SL_2\in \DD_c(\mls Y, \Lambda )$ denote $p_1^*L_1\Otimes p_2^*L_2$.

\begin{lem}\label{9.5} There is a natural isomorphism $\Omega _{\mls Y}\simeq \Omega _{\mls Y_1}\Otimes _S\Omega _{\mls Y_2}$ in $\DD_c(\mls Y, \Lambda )$.
\end{lem}
\begin{proof}
This is reduced to \cite[5.7.1]{Las-Ols05} by the same argument proving \ref{def-Omega-adic} using the gluing lemma.
\end{proof}

\begin{lem}\label{9.6b} For $L_i\in \DD_c^-(\mls Y_i)$ ($i=1,2$) there is a canonical isomorphism
\begin{equation}\label{dualarrow}
D_{\mls Y_1}(L_1)\Otimes _SD_{\mls Y_2}(L_2)\rightarrow D_{\mls Y}(L_1\Otimes _SL_2).
\end{equation}
\end{lem}
\begin{proof} Note first that there is a canonical map
\begin{equation}\label{pullmap}
p_i^*D_{\mls Y_i}(L_i)\rightarrow \Rghom (p_i^*L_i, p_i^*\Omega _{\mls Y_i}).
\end{equation}
Indeed by adjunction giving such a morphism is equivalent to giving a morphism
\begin{equation*}
D_{\mls Y_i}(L_i)\rightarrow Rp_{i*}\Rghom (p_i^*L_i, p_i^*\Omega _{\mls Y_i}),
\end{equation*}
and this in turn is by \ref{8.3} equivalent to giving a morphism
\begin{equation*}
D_{\mls Y_i}(L_i)\rightarrow \Rghom (L_i, Rp_{i*}p_i^*\Omega _{\mls Y_i}).
\end{equation*}
We therefore obtain the map \ref{pullmap} from the adjunction morphism $\Omega _{\mls Y_i}\rightarrow Rp_{i*}p_i^*\Omega _{\mls Y_i}$.

Combining this with \ref{9.4} we obtain a morphism
\begin{equation*}
D_{\mls Y_1}(L_1)\Otimes _SD_{\mls Y_2}(L_2)\rightarrow \Rghom (L_1\Otimes _SL_2, \Omega _{\mls Y_1}\Otimes _S\Omega _{\mls Y_2}),
\end{equation*}
which by \ref{9.5} defines the morphism \ref{pullmap}.

To see that this morphism is an isomorphism, note that by the definition this morphism is given by the natural map
\begin{equation*}
\Rhom _{\Lambda _\bullet }(\hat L_1, \Omega _{\mls Y_1})\Otimes _S\Rhom _{\Lambda _\bullet }(\hat L_2, \Omega _{\mls Y_2})\rightarrow \Rhom _{\Lambda _{\bullet }}(\hat L_1\Otimes _S\hat L_2, \Omega _{\mls Y})
\end{equation*}
in the topos $\mls Y^{\mathbb{N}}$.   That it is an isomorphism therefore follows from \cite[5.7.5]{Las-Ols05}.
\end{proof}

Let $f_i:\mls X_i\rightarrow \mls Y_i$ be morphisms of nice $S$--stacks, set $\mls X:= \mls X_1\times \mls X_2$, and let $f:= f_1\times f_2:\mls X_1\times\mls X_2\rightarrow \mls Y_1\times \mls Y_2$. Let $L_i\in \DD_c^-(\mls X_i, \Lambda)$ ($i=1,2$).

\begin{thm} There is a canonical isomorphism in $\DD_c(\mls Y, \Lambda )$
\begin{equation}
Rf_!(L_1\Otimes _SL_2)\rightarrow Rf_{1!}(L_1)\Otimes _SRf_{2!}(L_2).
\end{equation}
\end{thm}
\begin{proof}
As in \cite[proof of 5.7.5]{Las-Ols05} we define the morphism as the composite
\begin{equation*}
\begin{CD}
Rf_!(L_1{\Otimes_S}L_2)@>\simeq >> D_{\mc Y}(f_*D_{\mc X}(L_1{\Otimes_S}L_2))\\
@>\simeq >> D_{\mc Y}(f_*(D_{\mc X_1}(L_1){\Otimes_S}D_{\mc X_2}(L_2)))\\
@>>>  D_{\mc Y}(f_{1*}D_{\mc X_1}(L_1){\Otimes_S}(f_{2*}D_{\mc X_2}(L_2)))\\
@>\simeq >> D_{\mc Y_1}(f_{1*}D_{\mc X_1}(L_1)){\Otimes_S}D_{\mc Y_2}(f_{2*}D_{\mc X_2}(L_2))\\
@>\simeq >> Rf_{1!}(L_1){\Otimes_S}Rf_{2!}(L_2).
\end{CD}
\end{equation*}
That this morphism is an isomorphism is reduced, as in the proof of \ref{9.6b}, to loc. cit.
\end{proof}

\section{Base change theorem}

\begin{thm}\label{basechangethm} Let
\begin{equation*}
\begin{CD}
\mls X'@>a>> \mls X\\
@Vf'VV @VVfV \\
\mls Y'@>b>> \mls Y
\end{CD}
\end{equation*}
be a cartesian square of nice $S$--stacks.  Then the two functors
\begin{equation*}
b^*Rf_!, Rf'_!a^*:\DD _c^-(\mls X, \Lambda )\rightarrow \DD_c^-(\mls Y', \Lambda )
\end{equation*}
are canonically isomorphic.
\end{thm}

The proof of \ref{basechangethm} follows essentially the same outline as the proof of \cite[5.5.6]{Las-Ols05}.

\begin{lem} For any $A, B, C\in \DD_c(\X, \Lambda )$ there is a canonical morphism
$$
A\Otimes \Rghom (B, C)\rightarrow \Rghom (\Rghom (A, B), C).
$$
\end{lem}
\begin{proof} This is shown by the same argument proving \cite[5.5.7]{Las-Ols05}.
\end{proof}
By \cite[5.4.4]{Las-Ols05}, there exists a commutative diagram
\begin{equation}\label{5.5.6.3}
\begin{CD}
Y_\bullet '@>j>> Y_\bullet \\
@VpVV @VVqV \\
\mc  Y'@>b >> \mc  Y,
\end{CD}
\end{equation}
where $p$ and $q$ are smooth hypercovers and $j$ is a closed
immersion.

Let $\mc  X'_{Y'_\bullet }$ denote the base change $\mc  X'\times
_{\mc  Y'}Y_{\bullet }'$ and $\mc  X_{Y_\bullet }$ the base change
$\mc  X\times _{\mc  Y}Y_\bullet $. Then there is a cartesian
diagram
\begin{equation*}
\begin{CD}
\mc  X'_{Y'_\bullet }@>i>> \mc  X_{Y_\bullet }\\
@Vg'VV @VVgV \\
Y'_\bullet @>j>> Y_\bullet ,
\end{CD}
\end{equation*}
where $i$ and $j$ are closed immersions.

As in \cite[section 5]{Las-Ols05} let $\omega _{\X ' _{Y'_\bullet }}$ (resp. $\omega
_{\X_{Y_\bullet }}$, $\omega _{Y'_\bullet }$, $\omega _{Y_\bullet }$)
denote the pullback of $\Omega _{\X'}$ (resp. $\Omega _\X$,
$\Omega _{\Y'}$, $\Omega _{\Y}$) to $\X'_{Y'_\bullet }$ (resp.
$\X_{Y_\bullet }$, $Y_\bullet '$, $Y_\bullet $), and let
$D_{\X'_{Y'_\bullet }}$ (resp. $D_{\X_{Y_\bullet }}$, $D_{Y'_\bullet }$,
$D_{Y_\bullet }$) denote the functor $\Rhom (-, \omega
_{X'_{Y'_\bullet }})$ (resp. $\Rhom (-, \omega _{\X_{Y_\bullet }})$,
$\Rhom (-, \omega _{Y'_\bullet })$, $\Rhom (-, \omega _{Y_\bullet
})$). Note that these functors are defined already on the level of the derivated category of projective systems, though they also pass to the categories $\DD_c$.

Let $\mls F$ denote the functor
$$
D_{Y'_\bullet }j^*D_{Y_\bullet }Rg_*D_{\mls X_{Y_\bullet }}i_*D_{\mls X'_{Y'_\bullet }}:\D_c(\mls X_{Y'_\bullet }^{\prime \mathbb{N}}, \Lambda _\bullet )\rightarrow \D_c(Y^{\prime \mathbb{N}}_\bullet , \Lambda _\bullet ).
$$
By the same argument proving \cite[5.5.8]{Las-Ols05} one sees that there is a canonical isomorphism $\mls F\simeq Rg'_*$ (define a morphism of functors as in the beginning of the proof of \cite[5.5.8]{Las-Ols05} and then to check that it is an isomorphism it suffices to consider each $\Lambda _n$ where the result is loc. cit.), and hence also an isomorphism $\textbf{F}\simeq Rg'_*:\DD_c(\mls X'_{Y'_\bullet }, \Lambda )\rightarrow \DD_c(Y'_{\bullet }, \Lambda )$. This isomorphism induces a morphism of functors
\begin{equation}\label{firstmap}
\begin{matrix}
j^*D_{Y_\bullet }Rg_*D_{\X_{Y_\bullet }} & \rightarrow & j^*D_{Y_\bullet }Rg_*D_{\X_{Y_\bullet }}i_*i^* \ \ \  (\text{id}\rightarrow i_*i^*)\\
& \simeq & D_{Y'_\bullet }D_{Y'_\bullet }j^*D_{Y_\bullet }Rg_*D_{\X_{Y_\bullet }}i_*D_{\X'_{Y'_\bullet }}D_{\X'_{Y'_\bullet }}i^*\ \ \ \  (\ref{dual-adic})\\
& \simeq & D_{Y'_\bullet }\textbf{F}D_{\X'_{Y'_\bullet }}i^*\ \ \ \  (\text{definition})\\
& \simeq & D_{Y'_\bullet }Rg'_*D_{\X'_{Y'_\bullet }}i^*.
\end{matrix}
\end{equation}
If $\epsilon :\mls X_{Y_\bullet }\rightarrow \mls X$ and $\epsilon ':\mls X'_{Y_\bullet }\rightarrow \mls X'$ are the projections, we then obtain a morphism
\begin{equation}
\begin{matrix}
b ^*Rf_! & \simeq & b ^*Rq_*D_{Y_\bullet }Rg_*D_{\X_{Y_\bullet }}\epsilon ^*\ \ \ \ (\text{cohomological descent})\\
& \rightarrow & Rp_*j^*D_{Y_\bullet }Rg_*D_{\X_{Y_\bullet }}\epsilon ^*\ \ \ \  (\text{base change morphism})\\
&\rightarrow & Rp_*D_{Y'_\bullet }Rg'_*D_{\X'_{Y'_\bullet }}i^*\epsilon ^*\ \ \ \ (\ref{firstmap})\\
& \simeq & Rp_*D_{Y'_\bullet }Rg'_*D_{\X'_{Y'_\bullet }}\epsilon ^{\prime *}a^*\ \ \ \ (i^*\epsilon ^* = \epsilon ^{\prime *}a^*)\\
& \simeq & Rf'_!a^*\ \ \ \  (\text{cohomological descent}).
\end{matrix}
\end{equation}
which we call the \emph{base change morphism}. That it is an isomorphism is shown as in the proof of \cite[5.5.6]{Las-Ols05} by reduction to the case of schemes.
This completes the proof of \ref{basechangethm}. \qed

As in \cite{Las-Ols05} for some special classes of morphisms one can describe the base change arrow more explicitly.  The proofs that the following alternate definitions coincide with the one in \ref{basechangethm} proceed as in \cite{Las-Ols05} so we omit them.

By \ref{dual-adic} to prove the formula $b^*Rf_! = Rf'_!a^*$ it suffices to prove the dual version $b^!Rf_* = Rf'_*a^!$ which can do directly in the following cases.

\subsection{Smooth base change}

By \ref{9.3} the formula $b^!Rf_* = Rf'_*a^!$ is equivelent to the formula $b^*Rf_* = Rf'_*a^*$. We can therefore take the base change morphism $b^*Rf_*\rightarrow Rf'_*a^*$ (note that the construction of this arrow uses only adjunction for $(b^*, Rb_*)$ and $(a^*, Ra_*)$). To prove that this map is an isomorphism, note that it suffices to verify that it is an isomorphism locally in the topos $\mls Y^{\prime \mathbb{N}}$ where it follows from the case of finite coefficients \cite[5.1]{Las-Ols05}.

\subsection{Base change by a universal homeomorphism}

By \ref{9.6} in this case $b^! = b^*$ and $a^! = a^*$.  We then again take the base change arrow $b^*Rf_*\rightarrow Rf'_*a^*$ which as in the case of a smooth base change is an isomorphism by reduction to the case of finite coefficients \cite[5.4]{Las-Ols05}.

\subsection{Base change by  an immersion}

In this case one can define the base change arrow using the projection formula \ref{projection} as in \cite[5.3]{Las-Ols05}.

Note first of all that by shrinking on $\mls Y$ it suffices to consider the case of a closed immersion.
Let $A\in \DD_c(\mls X, \Lambda )$.  Since $b$ is a closed immersion, we have $b^*Rb_* = \text{id}$.  By the projection formulas for $b$ and $f$ we have
$$
Rb_*b^*Rf_!A = Rb_*(\Lambda )\Otimes Rf_!A = Rf_!(A\Otimes f^*Rb_*\Lambda ).
$$
Now clearly $f^*b_* = a_*f^{\prime *}$. We therefore have
\begin{equation*}
\begin{matrix}
Rf_!(A\Otimes f^*Rb_*\Lambda ) & \simeq & Rf_!(A\Otimes Ra_*f^{\prime *}\Lambda )\\
& \simeq & Rf_!a_*(a^*A\Otimes f^{\prime *}\Lambda ) \\
& \simeq & b_*Rf^{\prime }_!(a^*A)
\end{matrix}
\end{equation*}
Applying $b^*$ we obtain the base change isomorphism.

\end{document}